\documentclass[a4paper]{article}
\usepackage[english]{babel}
\usepackage{graphicx} 
\usepackage[T1]{fontenc}
\usepackage{csquotes}
\usepackage[backend=biber,
    giveninits=true,
    style=numeric-comp,
    maxbibnames=10,
    sorting=nyt,
    date=year,
    eprint=false,
    url=false
]{biblatex} 
\usepackage{amssymb}
\usepackage{amsxtra}
\usepackage{preamble}
\usepackage{todonotes}
\usepackage{tabularx}
\usepackage{easytable}
\usepackage{mathtools}
\usepackage[affil-it]{authblk}
\usepackage[a4paper,
            bindingoffset=0.2in,
            left=1in,
            right=1in,
            top=1in,
            bottom=1in,
            footskip=.25in]{geometry}

\usepackage{hyperref}
\hypersetup{hypertexnames=false}

\title{Abelian approximations for bicrossed products}
\author{Arthur Massar%
    \thanks{\texttt{arthur.massar@uclouvain.be}}}
\affil{Institut de Recherche en Math\'ematique et Physique, Universit\'e Catholique de Louvain, Chemin du Cyclotron 2, 1348 Louvain-La-Neuve, Belgique}
\date{\today}

\begin{document}

\maketitle

\begin{abstract}
    We propose a general framework to interpret the cocycle bicrossed product of a matched pair \(G \bowtie H = D\) of \lcgs as a deformation quantization of an extension \(V \to \Gcal \to G\) with \(V\) abelian.
    This construction requires the existence of a \(G\)-equivariant measure class isomorphism \(\eta_H : H \to \hat V\), which we call an \emph{abelian approximation} of \(H\).
    In the smooth case, the extension \(\Gcal\) naturally becomes equipped with a Lie--Poisson structure.
    When both \(G\) and \(H\) admit abelian approximations, the duality of the bicrossed products provides a unitary operator \(\Fcal : \Lrm^2(\Gcal) \to \Lrm^2(\Gcal^*)\), the \emph{Poisson--Fourier transform}, which implements the quantum duality principle of Drinfel'd.
    We discuss the construction of abelian approximations from Lie bialgebras and present several classes of examples. 
\end{abstract}

\section{Introduction}

Duality has always been a major motivation in the development of operator algebraic approaches to quantum groups since its very beginning, as one of its first objectives was the extension of Pontryagin duality to non-abelian locally compact groups \cite{enockKacAlgebrasDuality1992,vainermanNonunimodularRinggroupsHopfvon1974}.
One famous such approach is the theory of locally compact quantum groups of Kustermans and Vaes \cite{kustermansLCQG2000,kustermansvNLCQG2003}, which is the framework we will use in this paper.
Any \lcqg $(M, \Delta)$ admits a dual \lcqg $(\hat M, \hat \Delta)$, and there is a canonical isomorphism $(\hathat{M}, \hathat{\Delta}) \cong (M, \Delta)$.
We will also use $(M, \Delta)\sphat\,$ or $M\sphat\,$ to denote the dual $(\hat M, \hat \Delta)$ when this results in less cluttered notation.
A \lcg $\Gcal$ with left Haar measure $\dd g$ is represented by $M = \Lrm^{\infty}(\Gcal)$ and
\begin{align*}
    &\Delta_{\Gcal} : \Lrm^{\infty}(\Gcal) \to \Lrm^{\infty}(\Gcal) \htensor \Lrm^{\infty}(\Gcal) \cong \Lrm^{\infty}(\Gcal \times \Gcal) \\
    &\Delta_\Gcal(f)(g_1, g_2) = f(g_1 g_2).
\end{align*}
Its dual is modelled by the von Neumann algebra $\hat M = W^*(\Gcal) = \overline{\mathrm{span}\{\lambda_g : g \in \Gcal\}^{\text{WOT}}}$, equipped with the
\begin{align*}
    &\hat \Delta_{\Gcal} : W^*(\Gcal) \to W^*(\Gcal) \htensor W^*(\Gcal) \\
    &\hat \Delta_{\Gcal}(\lambda_g) = \lambda_g \tensor \lambda_g.
\end{align*}
Here $\lambda : \Gcal \to \Bcal(\Lrm^2(\Gcal))$ denotes the left regular representation.
For $\Gcal$ is abelian, the Fourier transform
\[
    \Fcal_{\Gcal} : \Lrm^2(\Gcal) \to \Lrm^2(\hat \Gcal)
\]
induces an isomorphism $\Fcal_G : (W^*(\Gcal), \hat \Delta_{\Gcal}) \cong (\Lrm^{\infty}(\hat \Gcal), \Delta_{\hat \Gcal})$, which recovers Pontryagin duality.

Our main goal in this paper will be to construct analogous Fourier transforms for some classes of \emph{non-abelian} groups.
For a non-abelian \lcg $\Gcal$ there clearly can not exist a \lcg $\Gcal'$ and a $*$-isomorphism $W^*(\Gcal) \to \Lrm^{\infty}(\Gcal')$, since $\Lrm^{\infty}(\Gcal')$ is commutative while $W^*(\Gcal)$ is not.
Nevertheless this becomes plausible again after quantization, where we deform $\Lrm^{\infty}(\Gcal')$ to a continuous family $\{\Lrm^{\infty}(\Gcal')_{\hbar}\}_{\hbar \in \R}$ of \lcqgs with $\Lrm^{\infty}(\Gcal')_0 \cong \Lrm^{\infty}(\Gcal')$ and $\Lrm^{\infty}(\Gcal')_{\hbar}$ non-commutative for $\hbar \neq 0$.
This is essentially a reformulation of the \qdp of Drinfel'd \cite{drinfeldQuantumGroups1988}:
Consider a Poisson--Lie group $(\Gcal, \{\cdot, \cdot\})$ with Poisson dual $(\Gcal^*, \{\cdot, \cdot\}_{*})$, and suppose we can quantize $(\Gcal, \{\cdot, \cdot\})$ to a family $\{\Lrm^{\infty}(\Gcal)_{\hbar}\}_{\hbar \in \R}$.
Then there should exist a corresponding quantization $\{\Lrm^{\infty}(\Gcal^*)_{\hbar}\}_{\hbar\in \R}$ of $(\Gcal^*, \{\cdot, \cdot\}_{*})$ and a family of isomorphisms of \lcqgs
\begin{align*}
    \Fcal_{\hbar} : \Lrm^{\infty}(\Gcal)_{\hbar}\sphat\, \xrightarrow{\cong} \Lrm^{\infty}(\Gcal^*)_{\hbar},
    && \forall \hbar \in \R \setminus \{0\}.
\end{align*}
Here $\Lrm^{\infty}(\Gcal)_{\hbar}\sphat\,$ denotes the \lcqg dual to $\Lrm^{\infty}(\Gcal)_{\hbar}$.
We call $\{\Fcal_{\hbar}\}_{\hbar}$ the \emph{Poisson--Fourier transform} between $\{\Lrm^{\infty}(\Gcal)_{\hbar}\}_{\hbar}$ and $\{\Lrm^{\infty}(\Gcal^*)_{\hbar}\}_{\hbar}$.
Drinfel'd provided a first construction of this duality for quantized Lie bialgebras and quantized polynomial function algebras, which was successfully exploited by Gavarini to prove general quantum duality theorems \cite{gavariniQuantumDualityPrinciple2002,gavariniGlobalQuantumDuality2007}.
We also recommend \cite{semenov-tyan-shanskiiPoissonLieGroupsQuantum1992} for an inspiring discussion.

So far the \qdp has not gotten much consideration in the operator algebraic framework, although it has been used to construct examples of \lcqgs \cite{rieffelSolvableQuantumGroups1989,SZ90QuantumHeisenberg,vandaeleQuantumDeformationHeisenberg1991} and a few instances of this duality have been exhibited, such as for the quantum $E(2)$ group \cite{woronowiczQuantumE2Group1991,vandaeleDualityQuantumE21996} and the quantum $az+b$ group \cite{woronowiczQuantumAzGroup2001} and the quantum $ax+b$ group \cite{woronowiczQuantumAxGroup2002,puszNewQuantumDeformation2005}.
However the \qdp as stated above poses a few problems which do not arise at the Lie bialgebra level due to the ambiguity of Poisson--Lie duality: $\Gcal^*$ is only determined by its Lie algebra and is therefore not unique.
Given a quantization $\{\Lrm^{\infty}(\Gcal)_{\hbar}\}$ of $(\Gcal, \{\cdot, \cdot\})$ we expect there to be a single choice of Poisson dual $\Gcal^*$ for which the \qdp makes sense, since the duality of \lcqgs is unique.
We call this specific choice the \emph{preferred} Poisson dual, and denote it $\Gcal^*_{\mathrm{pref}}$.
The question then is, how can we characterise $\Gcal^*_{\mathrm{pref}}$ among all the Poisson dual groups?
During an interesting discussion with Kenny de Commer, we also asked ourselves whether it could be possible to find two non-isomorphic quantizations $\{\Lrm^{\infty}(\Gcal)_{\hbar}\}_{\hbar}$ and $\{\Lrm^{\infty}(\Gcal)_{\hbar}'\}_{\hbar}$ of the same Poisson--Lie group $(\Gcal, \{\cdot,\cdot\})$ for which the preferred duals $\Gcal^*_{\mathrm{pref}}$ and $\Gcal^{*\prime}_{\mathrm{pref}}$ are \emph{not} isomorphic.
It is common in the Poisson--Lie group literature to choose the connected and simply connected dual, but as has been observed several times it can be necessary to consider non-connected groups for the purposes of quantization.
Notable examples include the quantum $\SU(1,1) \rtimes \Z/2$ of Korogodsky \cite{korogodskyQuantumSU[1-1]1994} and Koelink--Kustermans \cite{koelinkLocallyCompactQuantum2003}, as well as the quantum $ax+b$ group of Baaj and Skandalis \cite{skandalisDualityLocallyCompact1991}.
Connexity is also the key reason why the dual $2$-cocycles constructed by Bieliavsky and Gayral on negatively curved K\"alerian Lie groups in \cite{bieliavskyDeformationQuantizationActions2015} turned out to be coisometric but not unitary \cite[Appendix B. Erratum]{bieliavskyDeformationsCalgebrasActions2016}.
Similarly we will sometimes need to consider non-connected groups to make sense of the \qdp.

In this paper we will introduce a general framework to understand bicrossed products of \lcgs under the light of the \qdp.
The bicrossed product of Hopf algebras was introduced by Majid \cite{majidPhysicsAlgebraistsNoncommutative1990}, generalising earlier work of Kac \cite{kacExtensionsGroupsRing1968}, Singer \cite{singerExtensionTheoryConnected1972}, and Takeuchi \cite{takeuchiMatchedPairsGroups1981}.
It has since become a key tool in the theory of Hopf algebras.
Already in the nineties, examples of this construction were studied at the analytical level using the language of Hopf--von Neumann algebras \cite{majidHopfvonNeumannAlgebra1991}.
Many algebraic examples correspond to quantizations of classical groups or Lie algebras, such as the celebrated $\kappa$-Poincaré algebra \cite{majidBicrossedKappaPoincare1994}.
Unfortunately in the operator algebraic framework it is often unclear how to understand a bicrossed product as a quantization in the sense of deformation quantization, and answering this question will be one of our major objectives.
The ideas presented in this paper came to be through the study of the quantum $ax+b$ group of Baaj and Skandalis \cite{skandalisDualityLocallyCompact1991}, which is the bicrossed product of two copies of $\R^{\times}$ inside $\R \rtimes \R^{\times}$.
At first glance it is not obvious why this should correspond to a quantization of the $ax+b$ group.
 Vaes and Vainerman provided a first clue by analysing the tangent Hopf algebra of the bicrossed product and showing it forms a deformation of the Lie algebra of the $ax+b$ group \cite[Section 5.3]{vaesExtensionsLocallyCompact2003}.
 In 2013, Stachura constructed a unitary dual $2$-cocycle on the $ax+b$ group and conjectured the resulting quantization to be isomorphic to the bicrossed product of Baaj and Skandalis \cite{stachuraQuantum$axB$2013}.
 This question was answered positively by Bieliavsky, Gayral, Neshveyev, and Tuset \cite{bieliavskyQuantizationSubgroupsAffine2021}.
 The key idea was to use a Fourier transform on one of the components of the bicrossed product to compare its structure with the cocycle twist of $ax+b$.
 We will expand on this strategy and generalise it to any cocycle bicrossed products, provided one can construct suitable maps called \emph{abelian approximations}.
This provides both an approach to interpret bicrossed products as quantizations of classical groups in the sense of deformation quantization, and a strategy to obtain new examples of (cocycle) bicrossed products.

After a brief recollection of the bicrossed product construction in \zcref{section: bicrossed product}, we turn our attention to the semi-classical theory to understand which Poisson--Lie groups we may hope to quantize through these methods.
The main ingredient in the bicrossed product construction is a matched pair of Lie groups (or double Lie group) \(D = G \bowtie H\).
That is, \(D\) is a Lie group and \(G, H \leq D\) are Lie subgroups such that \(G \cap H = \{e_D\}\) and \(D \setminus GH\) has measure \(0\), see \zcref{def: matched pair}.
At the infinitesimal level, this corresponds to a \emph{double Lie algebra} \(\dfr = \gfr + \hfr\).
As has been observed before \cite[Proposition 8.3.3]{majidFoundationsQuantumGroup1995}, such double Lie algebras are in bijection with Lie bialgebras of the form \(\Gfr = \gfr \ltimes_{\alpha^{\flat}} \hfr^*_{\ab}, \Gfr^* = \hfr \ltimes_{\beta^{\flat}} \gfr^*_{\ab}\).
Here \(\gfr^*_{\ab}\) and \(\hfr^*_{\ab}\) represent the abelian Lie algebras on the vector spaces \(\gfr^*\), \(\hfr^*\), and \(\alpha : \gfr \to \gl(\hfr)\), \(\beta : \hfr \to \gl(\gfr)\) denote the infinitesimal dressing transformations defined by \([X, Y]_{\dfr} = \alpha_X Y - \beta_Y X\) for \(X \in \gfr\), \(Y \in \hfr\).
In particular if \(G\) is connected and simply connected, then the matched pair \(D = G \bowtie H\) induces a Lie--Poisson structure on the semi-direct product \(G \ltimes \hfr^*\).
In \zcref{section: lie bialgebras of translation type} we extend the above discussion to double Lie algebras with matched cocycles, and refer to the resulting Lie bialgebras (resp. Poisson--Lie groups) as Lie bialgebras (resp. Poisson--Lie groups) of \emph{translation type} since all of the non-commutativity of the dual is contained in the dual of the abelian ideal.
This terminology is also inspired by the name T-Poincaré groups of \cite{mercatiTMinkowski_I_2024}, see \zcref{ex: T-Poincare}.

\zcref[S]{section: abelian approximation} is the central part of the paper, where we investigate the quantization of extensions \(V \to \Gcal \to G\) with \(V\) abelian using a cocycle bicrossed product \(\hat G \bicrossl_{\alpha,\Ucal}^{\beta,\Vcal} H\).
For simplicity we restrict ourselves in the introduction to semi-direct products \(\Gcal = G \ltimes V\).
The following definition is the key to our construction.
\begin{definition}[{\zcref[S]{def:abelianApprox}}]
    Let $G \bowtie H = D$ be a matched pair of \lcgs with dressing transformations $\alpha : G \times H \to H, \beta : H \times G \to G$.
    That is, $\alpha$ and $\beta$ are determined almost-everywhere by the equation $g h^{-1} = \alpha_g(h)^{-1} \beta_{h}(g)$ in $D$.
    A \emph{measurable abelian approximation of $H$} consists of a \lcab $V$ and a measure class isomorphism $\eta_H : H \to \hat V$ such that the Borel action
    \begin{align*}
        &\hat \rho : G \times \hat V \to \hat V \\
        &\hat \rho_g(\xi) = \eta_H(\alpha_g(\eta_H^{-1}(\xi)))
    \end{align*}
    extends to a left action by continuous automorphisms.
    The \emph{induced action} is the right action $\rho : V \times G \to V$ dual to $\hat \rho$, in the sense that $\ppairing{\rho(v)_g}{\xi} = \ppairing{v}{\hat \rho_g(\xi)}$.
    Here, $\ppairing{\cdot}{\cdot} : V \times \hat V \to U(1)$ denotes the Pontryagin pairing.
    The \emph{classical limit} associated to $\eta_H$ is the semi-direct product $\Gcal = G \ltimes_{\rho} V$.
\end{definition}
If $\eta_H : H \to \hat V$ is an abelian approximation, we then have a chain of isomorphisms of von  Neumann algebra
\[
    W^*(G) \ltimes_{\alpha} \Lrm^{\infty}(H) \xrightarrow{\id \tensor {\eta_{H}}_*} W^*(G) \ltimes_{\hat \rho} \Lrm^{\infty}(\hat V) \xrightarrow{\id \tensor \Fcal_V^*} W^*(\Gcal).
\]
We refer to the composition \(\Fcal_H := \id \tensor \Fcal_V^* {\eta_H}_* : W^*(G) \ltimes_{\alpha} \Lrm^{\infty}(H) \to W^*(\Gcal)\) as the \emph{Poisson--Fourier \(\frac{1}{2}\)-transform}.
We then transport the coproduct $\Delta^{\beta}$ on $W^*(G) \ltimes_{\alpha} \Lrm^{\infty}(H)$ through \(\Fcal_H\) to a new coproduct $\hat \Delta^{\beta}_{\Gcal}$ on $W^*(\Gcal)$, and interpret the corresponding \lcqg $(W^*(\Gcal), \hat \Delta^{\beta}_{\Gcal})$ as a quantization of $\Gcal$.
The coproduct \(\hat \Delta^{\beta}_{\Gcal}\) induces a associative product \(\star^{\beta}\) on the Fourier algebra \(\Acal(\Gcal)\) since \(\Acal(\Gcal)\) is the predual of \(W^*(\Gcal)\).
To really say that \((W^*(\Gcal), \hat \Delta^{\beta}_{\Gcal})\) is a quantization of \((W^*(\Gcal), \hat \Delta_{\Gcal})\), one should include a deformation parameter and perform the semi-classical analysis of the resulting family.
This is often possible in practice but a general consideration of this question is beyond the scope of this paper, see \zcref{rk: semi-classical limit}.

Now suppose we are given a matched pair $G \bowtie H = D$ together with two abelian approximations $\eta_H : H \to \hat V$ and \(\eta_G : G \to \hat W\) of $H$ and \(G\) respectively, with classical limits $\Gcal = G \ltimes_{\rho} V$ and \(\Hcal = H \ltimes_{\sigma} W\).
In this situation we say that $(\eta_G, \eta_H)$ are \emph{matched abelian approximations}, see \zcref{def:matchedAbApprox}.
We write $\hat \Delta^{\alpha}_{\Hcal}$ for the new coproduct on $W^*(\Hcal)$ constructed from $\eta_G$ using the procedure above.
Combining $\Fcal_H$ and $\Fcal_G$ with the flip map $\Sigma : \Lrm^2(G) \tensor \Lrm^2(H) \to \Lrm^2(H) \tensor \Lrm^2(G)$, we obtain a unitary operator
\[
    \Fcal_{\Gcal} : \Lrm^2(\Gcal) \xrightarrow{\Fcal_H^*} \Lrm^2(G) \tensor \Lrm^2(H) \xrightarrow{\Sigma} \Lrm^2(H) \tensor \Lrm^2(G) \xrightarrow{\Fcal_G} \Lrm^2(\Hcal),
\]
which induces an isomorphism of \lcqgs
\begin{align*}
    \Fcal_{\Gcal} : (W^*(\Gcal), \hat \Delta^{\beta}_{\Gcal})\sphat\; \cong (W^*(\Hcal), \hat \Delta^{\alpha}_{\Hcal})
\end{align*}
between the dual of the quantization $(W^*(\Gcal),\hat \Delta^{\beta}_{\Gcal})$ of $\Gcal$ and the quantization $(W^*(\Hcal), \hat\Delta^{\alpha}_{\Hcal})$ of $\Hcal$.
We call $\Fcal_{\Gcal}$ the \emph{Poisson--Fourier transform} between $\Gcal$ and $\Hcal$, and interpret the above isomorphism as implementing the quantum duality principle of Drinfel'd.

We conclude the theoretical part of this article with \zcref{section: abelian approximation from coboundary structure} where we present a strategy to produce abelian approximations from coboundary structures, which is inspired by the following example.
Consider a connected semi-simple real Lie group \(D\) with its Iwasawa and Cartan decompositions \(D = K \bowtie AN = KP\) and write \(\dfr = \kfr + (\afr + \nfr) = \kfr + \pfr\) for the Lie algebraic decompositions.
We define maps \(\Lrm^K_P : D \to K\) and \(\Rrm^P : D \to P\) by
\begin{align*}
    d = \Lrm^K_P(d) \Rrm^P(d) && \text{for all } d \in D.
\end{align*}
Since the restricted exponential map \(\exp|_{\pfr} : \pfr \to P\) is a \(K\)-equivariant diffeomorphism, the composition \(\eta := \exp|_{\pfr}^{-1} \circ \Rrm^P : AN \to \pfr\) is an abelian approximation for \(D = K \bowtie AN\).
The key properties which make the above construction possible are
\begin{enumerate}
    \item \(\Ad_g \pfr = \pfr\) for all \(g \in K\);
    \item \(\exp_{\pfr} : \pfr \to P\) is a \(K\)-equivariant diffeomorphism;
    \item \(D = K P\) and \(K \cap P = \{e_D\}\).
\end{enumerate}
An \emph{exponential reductive complement} for a double Lie group \(D = K \bowtie H\) is then the data of a linear subspace \(\pfr \subseteq \dfr\) and a submanifold \(P \subseteq D\) which satisfy the above axioms, see \zcref{def: reductive complement}.
In \zcref{prop: abelian approximation from lie--reductive doubles} we show that for any double Lie group with exponential reductive complement \(D = K \bowtie H = K P\), the map \(\eta_H := \exp|_{\pfr}^{-1} \circ \Rrm^P : H \to \pfr\) is an abelian approximation.
This compels us to construct such subspaces \(\pfr \subseteq \dfr\) as a first step towards constructing abelian approximations.
A careful analysis of the classical double of coboundary Lie bialgebras yields the following statement.
\begin{proposition}[{\zcref[S]{cor: reductive complement from coboundary of translation type}}]
    Let \((\Gfr = \gfr \ltimes_{\alpha} \hfr^*_{\ab}, r)\) be a coboundary Lie bialgebra of translation with dual \(\Gfr^* = \hfr \ltimes_{\beta} \gfr^*_{\ab}\) and write \(\dfr = \gfr + \hfr\) for the associated double Lie algebra.
    Then \(\dfr\) decomposes as \(\dfr = \gfr + \pfr_+ = \gfr + \pfr_-\) for linear subspaces \(\pfr_+, \pfr_- \subseteq \dfr\) such that \([\gfr, \pfr_{\pm}]_{\dfr} \subseteq \pfr_{\pm}\), and \([\pfr_+, \pfr_-]_{\dfr} \subseteq \gfr\).
    Moreover,
    \begin{enumerate}
        \item If \(r\) is skew-symmetric, then \(\pfr_+ = \pfr_- =: \pfr\) and \([\pfr, \pfr]_{\dfr} \subseteq \gfr\);
        \item If \(r\) is quasi-triangular, then \(\pfr_+, \pfr_-\) are ideals and \([\pfr_+, \pfr_-]_{\dfr} = 0\);
        \item If \(r\) is triangular, then \(\pfr_+ = \pfr_- =: \pfr\) is an abelian ideal and \(\dfr \cong \gfr \ltimes_{\alpha} \hfr_{\ab}\).
    \end{enumerate}
\end{proposition}
As a consequence we get a general quantization scheme for triangular Lie bialgebras of translation type.
\begin{corollary}[{\zcref[S]{cor: abelian approximation for triangular TT}}]
    Consider a Lie algebra \(\Gfr = \gfr \ltimes_{\alpha} V\) and a triangular \(r\)-matrix of translation type \(r = r_{\gfr,V} + r_{V,\gfr} \in \Gfr \wedge \Gfr\), where \(r_{\gfr,V} \in \gfr \tensor V\), \(r_{V,\gfr} \in V \tensor \gfr\).
    Then \(\hfr = \{f + r^{\sharp}(f) : f \in V^*\} \subseteq \gfr \ltimes_{\alpha^{\flat}} V^*\) is a Lie subalgebra.
    If the double Lie algebra \(\gfr \ltimes_{\alpha^{\flat}}V^* = \gfr + \hfr\) integrates to a double Lie group \(G \ltimes V^* = G \bowtie H\), then the map
    \begin{align*}
        &R^{V^*} : H \to V^* \\
        &R^{V^*}(h) := f && \text{for all } h = gf \in H \subseteq G \ltimes V^*\!, g \in G, f \in V^*
    \end{align*}
    is an abelian approximation.
\end{corollary}
An analogous statement holds in the quasi-triangular case when \(\pfr_+\) is the Lie algebra of an exponential Lie group as explained in \zcref{cor: abelian approximation for quasi-triangular TT}.

\zcref[S]{sec:examples} is dedicated to examples of the Poisson--Fourier transform, including the quantum \(ax+b\) group of Baaj and Skandalis and the quantum Heisenberg groups of Szymczak and Zakrzewski.

\subsection{Notations}

Let $G$ be a locally compact group.
We always assume $G$ to be second countable.
In particular, the left Haar measure $\dd g$ is $\sigma$-finite.
We let $\Lrm^p(G) = \Lrm^p(G, \C), p \in [1, \infty]$ denote the corresponding function spaces.
The modular function $\delta_G$ of $G$ is defined by the equation
\begin{equation}
    \int_G f(g g_0) \dd g = \delta_G(g_0)^{-1} \int_G f(g) \dd g
\end{equation}
for $f \in C_c(G)$.

Let $V$ be a locally compact abelian group, $\rho : V \times G \to V : (v, g) \mapsto \rho(v)_g$ a right action of $G$ on $V$ by continuous automorphisms, and $\vphi : G \times G \to V$ a $2$-cocycle for $\rho$.
We let $\Gcal = G \ltimes_{\rho,\vphi} V$ denote the extension of $G$ by $V$ determined by $\rho$ and $\vphi$.
That is, its group law is given by
\[
    (g_1, v_1)(g_2, v_2) = (g_1 g_2, \rho(v_1)_{g_2} + v_2 + \vphi(g_1, g_2)).
\]
If $\alpha : G \times V \to V$ is a left action instead, we denote the semi-direct product by $V \rtimes_{\alpha} G$.

We write $\hat V$ for the Pontryagin dual of $V$, and let 
\[
    \ppairing{-}{-} : V \times \hat V \to U(1)
\]
denote the Pontryagin pairing.
We normalise the Haar measures on $V$ and $\hat V$ such that the Fourier transform
\begin{align*}
    &\Fcal_V : \Lrm^2(V) \to \Lrm^2(\hat V), \\
    &\Fcal_V \vphi(\hat v) = \int_V \ppairing{\hat v}{v} \vphi(v) \dd v
\end{align*}
is unitary.
We let $\hat \rho : G \times \hat V \to \hat V : (g, v) \mapsto \hat\rho_g(v)$ denote the dual \emph{left} action of $G$ on $\hat V$, determined by
\[
    \ppairing{\rho(v)_g}{\hat v} = \ppairing{v}{\hat\rho_g(\hat v)}
\]
for $v \in V, \hat v \in \hat V, g \in G$.
As we will see, it is more natural in the framework to switch between right and left actions.

Suppose that $(X, \mu)$ and $(Y, \nu)$ are $\sigma$-finite measured spaces.
By a \emph{measure class isomorphism} we mean a measurable map $f : X \to Y$ which admits a measurable inverse $g : Y \to X$.
That is, we have $g \circ f = \id_X$ and $f \circ g = \id_Y$ almost everywhere, and the measures $f_* \mu$ and $\nu$ on $Y$ are absolutely continuous with respect to each other.
In particular this entails that the map
\begin{align*}
    &f^* : \Lrm^p(Y, \nu) \to \Lrm^p(X, \mu), \\
    &(f^*\vphi)(x) = \vphi(f(x)) \left(\frac{\dd g_* \nu}{\dd \mu}(x)\right)^{1/p}
\end{align*}
is an isomorphism for every $p \in [1, \infty]$.
We let $f_* := g^* : \Lrm^p(X, \mu) \to \Lrm^p(Y, \nu)$ denote its inverse, and set $f^{-1} := g$.

We write \(V_{\ab}\) for the abelian Lie algebra on a vector space \(V\).
We use the following terminology for Lie bialgebras:
\begin{definition}
    Consider a Lie bialgebra \(\gfr\) with cocommutator \(\gamma\).
    A \emph{coboundary structure} for \((\gfr, \gamma)\) is an element \(r \in \gfr \tensor \gfr\) such that \(\delta r = \gamma\).
    Then \(r\) is \emph{quasi-triangular} if it satisfies the classical Yang--Baxter equation
    \begin{align*}
        [r_{12}, r_{13}] + [r_{12}, r_{23}] + [r_{13}, r_{23}] = 0.
    \end{align*}
    It is \emph{triangular} if it is quasi-triangular and skew-symmetric.
\end{definition}

\section{The bicrossed product construction}\label{section: bicrossed product}

We begin by recalling the bicrossed product construction with cocycles.
We mostly follow \cite[Section 4.2]{vaesExtensionsLocallyCompact2003}, but we use slightly different notation when it comes to the cocycles.

 \begin{definition}\label{def: matched pair}
    A \emph{matched pair $G \bowtie H = D$} of \lcgs consists of a triple $(G, H, D)$, where $D$ is a \lcg and $G, H \subseteq D$ are closed subgroups such that $G \cap H = \{e_D\}$ and $GH = \{g h : g \in G, h \in h\}$ has complement of measure $0$ in $D$.
 \end{definition}

 Throughout this section we fix a matched pair $G \bowtie H = D$ and set
 \begin{align*}
    &\mu : G \times H \to D : (g,h) \mapsto g h^{-1};
    &\mu_{\Sigma} : G \times H \to D : (g, h) \mapsto h^{-1} g.
 \end{align*}

 By \cite[Proposition~3.2]{baajNonSemiRegularQuantumGroups2003} we can normalise the Haar measures on $D$, $G$, and $H$ such that for all positive Borel function $f : D \to \R_{\geq 0}$, we have
\begin{align*}
    \int_D f(x) \dd x
    &= \int_{G \times H} f(\mu(g,h)) \delta_D(h)^{-1} \dd g \dd h \\
    &= \int_{G \times H} f(\mu_{\Sigma}(g,h)) \delta_D(g) \delta_G(g)^{-1} \delta_H(h)^{-1} \dd g \dd h.
\end{align*}
In particular, $\mu$ and $\mu_{\Sigma}$ are measure class isomorphisms.

 \begin{definition}
    The \emph{dressing transformations} of the matched pair $G \bowtie H = D$ consist of the Borel maps $\alpha : G \times H \to H$ and $\beta : H \times G \to G$ determined by
    \[
        \mu_{\Sigma}^{-1}(\mu(g,h)) = (\beta_h(g), \alpha_g(h))
    \]
    for almost every $g \in G, h \in H$.
\end{definition}
Equivalently, $\alpha$ and $\beta$ are defined almost everywhere by the equation
\[
    g h^{-1} = \alpha_g(h)^{-1} \beta_h(g).
\]
From the associativity of the group structure of $D$, we see that $\alpha$ and $\beta$ satisfy the following identities for almost every $g, g' \in G$, $h, h' \in H$:
\begin{align}\label{eq:dresTransIds}
\begin{split}
    \alpha_{g' g}(h) = \alpha_{g'}(\alpha_g(h));\qquad
    &\beta_h(g' g) = \beta_{\alpha_g(h)}(g') \beta_h(g); \\
    \beta_{h' h}(g) = \beta_{h'}(\beta_h(g));\qquad
    &\alpha_g(h' h) = \alpha_{\beta_h(g)}(h') \alpha_g(h).
\end{split}
\end{align}
Moreover, we have
\begin{align*}
    \alpha_g(e_H) = e_h;
    && \beta_{e_H}(g) = g;
    && \alpha_{e_G}(h) = h; 
    && \beta_h(e_G) = e_G.
\end{align*}

\begin{definition}
    Let $G \bowtie H = D$ be a matched pair of \lcgs with dressing transformations $(\alpha, \beta)$.
    We say that Borel maps $\Ucal : G \times G \times H \to U(1)$ and $\Vcal : H \times H \times G \to U(1)$ form \emph{matched cocycles} for $G \bowtie H = D$ if they satisfy the following equations almost everywhere
    \begin{align}
        \label{eq:cocycleU} \Ucal(g_1, g_2, \alpha_{g_3}(h)) \Ucal(g_1 g_2, g_3, h)
        &= \Ucal(g_1, g_2 g_3, h) \Ucal(g_2, g_3, h), \\
        \label{eq:cocycleV} \Vcal(h_1, h_2, \beta_{h_3}(g)) \Vcal(h_1 h_2, h_3, g)
        &= \Vcal(h_1, h_2 h_3, g) \Vcal(h_2, h_C, g),
    \end{align}
    and
    \begin{equation}\label{eq:matchedCocycles}
    \begin{split}
        \Vcal(h_1, h_2, g_1 g_2) \cj{\Ucal}(g_1, g_2, h_1 h_2)
        ={} & \Vcal(\alpha_{\beta_{h_2}(g_2)}(h_1, \alpha_{g_2}(h_2), g_1))
            \Vcal(h_1, h_2, g_2) \\
        &{}\cj{\Ucal}(\beta_{\alpha_{g_2}(h_2)}(g_1), \beta_{h_2}(g_2), h_1)
            \cj{\Ucal}(g_1, g_2, h_2).
    \end{split}
    \end{equation}
\end{definition}

The \zcref{eq:cocycleU,eq:cocycleV} simply mean that the adjoint maps $\underline{\Ucal} : G \times G \to \Lrm^{\infty}(H)$ and $\underline{\Vcal} : H \times H \to \Lrm^{\infty}(G)$ are unitary $2$-cocycles for the actions induced by $\alpha$ and $\beta$.
The compatibility between $\Ucal$ and $\Vcal$ is entirely contained in \zcref{eq:matchedCocycles}.

\begin{remark}
    Define $\tau : \Lrm^{\infty}(G) \htensor \Lrm^{\infty}(H) \to \Lrm^{\infty}(G) \htensor \Lrm^{\infty}(H)$ by
        $\tau(\vphi)(g,h) = \vphi(\beta_h(g), \alpha_g(h))$.
    Then $(\tau, \Ucal, \Vcal_{321})$ is a cocycle matching between $(\Lrm^{\infty}(G), \Delta_G)$ and $(\Lrm^{\infty}(H), \Delta_H)$ in the sense of \cite[Definition 2.1]{vaesExtensionsLocallyCompact2003}.
    Here $\Vcal_{321}(g,h_1,h_2) := \Vcal(h_2,h_1,g)$.
\end{remark}

\begin{remark}
    As we will see in \zcref{rk: matched pairs and Poisson structures}, a matched pair of Lie groups \(G \bowtie H = D\) induces Poisson--Lie structures on the semi-direct products \(\Gcal = G \ltimes \hfr^*\), \(\Hcal = H \ltimes \gfr^*\).
    The addition of matched cocycles \((\Ucal, \Vcal)\) correspond to general extensions \(\hfr^* \to \Gcal \to G\) and \(\gfr^* \to \Hcal \to H\).
    We refer to the Poisson--Lie groups which arise this way as Poisson--Lie groups \emph{of translation type}, see \zcref{def: lie bialgebra of translation type}.
    There are important examples of Poisson--Lie groups of translation type for which the associated cocycles are non-trivial.
    For example the T-Poincaré models of \zcref{ex: T-Poincare} with \(\theta^{\mu \nu} \neq 0\) should be quantized by bicrossed products with non-trivial cocycles.
    This is an important motivation for including matched cocycles in our discussion.
\end{remark}

We set
\begin{align*}
    &\alpha^* : \Lrm^{\infty}(H) \to \Lrm^{\infty}(G) \htensor \Lrm^{\infty}(H) \\
    &\alpha^*(\vphi)(g, h) := \vphi(\alpha_g(h))
\end{align*}
almost everywhere.
We are now ready to define the bicrossed product of $G \bowtie H = D$ with cocycles $(\Ucal, \Vcal)$.
See \cite[Section 2.2]{vaesExtensionsLocallyCompact2003} for the proof that this indeed does define a \lcqg.

\begin{definition}
    Suppose that $G \bowtie H = D$ is a matched pair of \lcgs with cocycles $(\Ucal, \Vcal)$ are dressing transformations $(\alpha, \beta)$.
    We write \(W_G\xi(g_1, g_2) = \xi(g_1, g_1^{-1} g_2)\) for all \(\xi \in \Lrm^2(G)\).
    The \emph{cocycle bicrossed product} $\hat G \bicrossl_{\alpha,\Ucal}^{\beta,\Vcal} H$ is the \lcqg determined by the following data:
    \begin{enumerate}
        \item The underlying von Neumann algebra is the cocycle crossed product $W^*(G) \ltimes_{\alpha,\Ucal} \Lrm^{\infty}(H)$.
            That is, it is the von Neumann subalgebra of $\Bcal(\Lrm^2(G)) \htensor \Lrm^{\infty}(H)$ generated by
            \begin{align*}
                \alpha^*(\Lrm^{\infty}(H))
                && \text{and}
                && \left\{(\omega \tensor \id \tensor \id)(W_G \tensor 1)\Ucal^* \mid \omega \in \Lrm^1(G)\right\};
            \end{align*}
        \item The multiplicative unitary $W_{\alpha,\Ucal}^{\beta,\Vcal} \in \Bcal(\Lrm^2(G \times H \times G \times H))$ is given by the formula
        \begin{align*}
            W_{\alpha,\Ucal}^{\beta,\Vcal} \xi(g_1, h_1, g_2, h_2)
            ={} & \Ucal(\beta_{\alpha_{g_1}(h_1)^{-1}h_2}(g_2), g_1, h_1) \\
            &{} \cj{\Vcal}(\alpha_{g_1}(h_1),\alpha_{g_1}(h_1)^{-1} h_2, g_2) \\
            &{} \xi(\beta_{\alpha_{g_1}(h_1)^{-1} h_2}(g_2) g_1, h_1, g_2, \alpha_{g_1}(h_1)^{-1} h_2);
        \end{align*}
        \item The coproduct is induced by $W_{\alpha,\Ucal}^{\beta,\Vcal}$ in the usual way:
        \[
            \Delta^{\beta,\Vcal}(x) = {W_{\alpha,\Ucal}^{\beta,\Vcal}}^* (1 \tensor x) W_{\alpha,\Ucal}^{\beta,\Vcal}
        \]
        for $x \in W^*(G) \ltimes_{\alpha,\Ucal} \Lrm^{\infty}(H)$.
    \end{enumerate}
    The left Haar weight is the dual weight on the cocycle crossed product $W^*(G) \ltimes_{\alpha,\Ucal} \Lrm^{\infty}(H)$, see \cite[Definition 1.13]{vaesExtensionsLocallyCompact2003}.
\end{definition}

We conclude this section by giving a description of the dual quantum group of $\hat G \bicrossl_{\alpha,\Ucal}^{\beta,\Vcal} H$.

\begin{proposition}[{\cite[Theorem 2.13]{vaesExtensionsLocallyCompact2003}}]\label{prop:dualityBicrossed}
    Suppose that $G \bowtie H = D$ is a matched pair with dressing transformations $(\alpha, \beta)$ and matched cocycles $(\Ucal, \Vcal)$.
    Then the flip operator $\Sigma : \Lrm^2(G) \tensor \Lrm^2(H) \to \Lrm^2(H) \tensor \Lrm^2(G)$ induces an isomorphism of \lcqgs
    \[
        \Sigma : \left(\hat G \bicrossl_{\alpha,\Ucal}^{\beta,\Vcal} H\right)\sphat\, \cong \hat H \bicrossl_{\beta,\Vcal}^{\alpha,\Ucal} G.
    \]
\end{proposition}

\section{Lie bialgebras of translation type}\label{section: lie bialgebras of translation type}
The infinitesimal analogue of a matched pair \(G \bowtie H = D\) of Lie groups is a \emph{double Lie algebra}.
That is, it is a vector space decomposition \(\dfr = \gfr + \hfr\), where \(\gfr, \hfr \leq \dfr\) are two Lie subalgebras.
We define the \emph{infinitesimal dressing transformations} \(\alpha : \gfr \times \hfr \to \hfr\) and \(\beta : \hfr \times \gfr \to \gfr\) by
\begin{align*}
    [X, Y] = \alpha_X Y - \beta_Y X && \text{for all } X \in \gfr, Y \in \hfr.
\end{align*}
Then \(\alpha\) and \(\beta\) are Lie algebra representations, and they satisfy the equations
\begin{align} \label{equations: Lie alg dressing equations}
    \begin{split}
    \alpha_X([Y_1,  Y_2]) &= \alpha_{\beta_{Y_2}(X)}(Y_1) - \alpha_{\beta_{Y_1}(X)}(Y_2) + [Y_1, \alpha_X(Y_2)] - [Y_2, \alpha_X(Y_1)]; \\
    \beta_Y([X_1, X_2]) &= \beta_{\alpha_{X_2}(Y)}(X_1) - \beta_{\alpha_{X_1}(Y)}(X_2) + [X_1, \beta_Y(X_2)] - [X_2, \beta_{Y}(X_1)].
    \end{split}
\end{align}
These are the infinitesimal analogue of the equations \zcref{eq:dresTransIds} characterising dressing transformations.
Moreover any pair of Lie algebra representations \(\alpha' : \gfr \to \gl(\hfr)\), \(\beta' : \hfr \to \gl(\gfr)\) which satisfy \zcref[noname]{equations: Lie alg dressing equations} give rise to a double Lie algebra \(\dfr = \gfr + \hfr\), and these two constructions are inverses to each other \cite{kosmann-schwarzbachPoissonLieGroupsComplete1988}.

Similarly, the infinitesimal analogue of matched cocycle \(\Ucal, \Vcal\) for the matched pair \(G \bowtie H = D\) is a pair of \(2\)-cocycles \(\vphi \in Z^2(\gfr, (\hfr^*, \alpha^{\flat}))\), \(\psi \in Z^2(\hfr, (\gfr^*, \beta^{\flat}))\) which satisfy
\begin{align}\label{eq: Lie alg pentagonal 2-cocycle}
    \begin{split}
    0 &= \vphi(\beta_{Y_1} X_1, X_2)(Y_2) + \vphi(X_1, \beta_{Y_1} X_2)(Y_2)
        - \vphi(\beta_{Y_2} X_1, X_2)(Y_1) - \vphi(X_1, \beta_{Y_2} X_2)(Y_1) \\
        &\ + \psi(\alpha_{X_1} Y_1, Y_2)(X_2) + \psi(Y_1, \alpha_{X_1} Y_2)(X_2)
        - \psi(\alpha_{X_2} Y_1, Y_2)(X_1) - \psi(Y_1, \alpha_{X_2} Y_2)(X_1) \\
        &\ + \vphi(X_1, X_2)([Y_1, Y_2]) + \psi(Y_1, Y_2)([X_1, X_2])
    \end{split}
\end{align}
for \(X_i \in \gfr, Y_i \in \hfr, i \in \{1,2,3\}\).
Again this equation is related to \zcref{eq:matchedCocycles}.
The data of \(\vphi, \psi\) is equivalent to the data of a \(3\)-cocycle \(\Theta \in Z^3(\dfr, \R)\) such that \(\Theta|_{\gfr} = 0\) and \(\Theta|_{\hfr} = 0\), as explained in \cite[Proposition 1]{chapovskii22Locally2004}.
The relation between \(\Theta\) and \((\vphi, \psi)\) is
\begin{align} \label{eq: 3 cocycle from matched cocycles}
    \Theta(X, X', Y) = \la \vphi(X, X'), Y \ra; &&
    \Theta(Y, Y', X) = \la \psi(Y, Y'), X \ra
\end{align}
for \(X, X' \in \gfr\), \(Y, Y' \in \hfr\).

In this section we want to understand the semi-classical structure associated to a double Lie algebra \(\dfr = \gfr + \hfr\) equipped with matched cocycles \(\vphi, \psi\).
The resulting Lie bialgebras will be called \emph{Lie bialgebras of translation type}.
We first define Lie bialgebras of translation type and show the equivalence with double Lie algebras equipped with matched cocycles in \zcref{th: double lie algs and lie bialgs of translation type}.
When the cocycles \(\vphi, \psi\) are trivial, this construction reduces to a previous result of Majid \cite[Proposition 8.3.3]{majidFoundationsQuantumGroup1995}.

\begin{definition}\label{def: lie bialgebra of translation type}
    Consider two Lie algebras \(\gfr, \hfr\) together with representations \(\alpha : \gfr \to \gl(\hfr)\), \(\beta : \hfr \to \gl(\gfr)\) and \(2\)-cocycles \(\vphi \in Z^2(\gfr, (\hfr^*, \alpha^{\flat}))\), \(\psi \in Z^2(\hfr, (\gfr^*, \beta^{\flat}))\).
    If the extensions \(\Gfr = \gfr \ltimes_{\alpha^{\flat},\vphi} \hfr^*_{\ab}\) and \(\Hfr = \hfr \ltimes_{\beta^{\flat},\psi} \gfr^*_{\ab}\) form Lie bialgebras which are dual to each other with respect to the canonical pairing
    \begin{align*}
        &\la -, - \ra : \Gfr \times \Hfr \to \R \\
        &\la X + \eta, Y + \xi \ra = \la X, \xi \ra + \la Y, \eta \ra && \text{for all } X \in \gfr, \eta \in \hfr^*, Y \in \hfr, \xi \in \gfr^*,
    \end{align*}
    then we say that \(\Gfr, \Hfr\) are \emph{Lie bialgebras of translation type}.
\end{definition}

Note that all of the compatibility equations in the previous definition are contained in the condition that \((\Gfr, \Hfr)\) forms a Lie bialgebra.
We first characterise Lie bialgebras of translation type as those whose double \(\Dfr\) is an extension \(\Dfr = \dfr \ltimes_{\ad^{\flat},\omega} \dfr^*_{\ab}\).

\begin{proposition}\label{prop: double of lie bialg of translation type}
    Consider a Lie bialgebra of translation type \(\Gfr = \gfr \ltimes_{\alpha^{\flat},\vphi} \hfr^*_{\ab}\) with dual \(\Hfr = \hfr \ltimes_{\beta^{\flat}, \psi} \gfr^*_{\ab}\).
    Then their double \(\Dfr\) decomposes as an extension \(\Dfr = \dfr \ltimes_{\ad^{\flat},\omega} \dfr^*_{\ab}\), where \(\dfr = \gfr + \hfr\) is a Lie algebra with bracket
    \begin{align*}
        [X_1 + Y_1, X_2 + Y_2]_{\dfr} = [X_1, X_2]_{\gfr} + \alpha_{X_1} Y_2 - \beta_{Y_2} X_1 + \beta_{Y_1} X_2 - \alpha_{X_2} Y_1 + [Y_1, Y_2]_{\hfr},
    \end{align*}
    and \(\omega \in Z^2(\dfr, \dfr^*)\) is the \(2\)-cocycle defined by
    \begin{align*}
        \omega(X_1 + Y_1, X_2 + Y_2) = \vphi(X_1, X_2) + \psi(Y_1, Y_2) &- \vphi(X_1, -)^*(Y_2) + \psi(Y_2, -)^*(X_1) \\
        &+ \vphi(X_2, -)^*(Y_1) - \psi(Y_1, -)^*(X_2),
    \end{align*}
    where \(X_1, X_2 \in \gfr\), \(Y_1, Y_2 \in \hfr\).
    Conversely if \((\Dfr = \dfr \ltimes_{\ad^{\flat}, \omega} \dfr^*_{\ab}, \Gfr, \Hfr)\) is a Manin triple with respect to the inner product
    \(\la D_1 + \phi_1, D_2 + \phi_2\ra_{\Dfr} = \la D_1, \phi_2 \ra + \la D_2, \phi_1 \ra\) for all \(D_1, D_2 \in \dfr\), \(\phi_1, \phi_2 \in \dfr^*\), then \(\Gfr\) and \(\Hfr\) are Lie bialgebras of translation type which are dual to each other.
\end{proposition}
\begin{proof}
    We first compute the coadjoint action of \(\Gfr\):
    \begin{align}\label{eq: coadjoint action of extension}
        \begin{split}
        &\la \ad^{\flat}_{X_0 + \eta_0}(Y + \xi), X_1 + \eta_1 \ra
        = - \la Y + \xi, [X_0 + \eta_0, X_1 + \eta_1]_{\Gfr} \ra \\
        &\qquad= - \la \xi, [X_0, X_1]_{\gfr} \ra - \la Y, \vphi(X_0, X_1) + \alpha^{\flat}_{X_0} \eta_1 - \alpha^{\flat}_{X_1} \eta_0 \ra \\
        &\qquad= \la \ad^{\flat}_{X_0} \xi, X_1 \ra - \la \vphi(X_0, -)^*(Y), X_1 \ra + \la \alpha_{X_0} Y, \eta_1 \ra + \la (\alpha^{\flat}_{-} \eta_0)^*(Y), X_1 \ra,
        \end{split}
    \end{align}
    where \(X_0, X_1 \in \gfr, \eta_0, \eta_1 \in \hfr^*, Y \in \hfr, \xi \in \gfr^*\).
    Using the analogous formulas for \(\Hfr\) yields
    \begin{align}\label{eq: bracket on Dfr}
        \begin{split}
        [X + \eta, Y + \xi]_{\Dfr} &= \alpha_X Y - \beta_Y X + \ad^{\flat}_X \xi - \ad^{\flat}_Y \eta + (\alpha^{\flat}_{-} \eta)^*(Y) - (\beta^{\flat}_{-} \xi)^*(X) \\
        &\qquad- \vphi(X, -)^*(Y) + \psi(Y, -)^*(X).
        \end{split}
    \end{align}
    In particular \(\dfr^* \subseteq \Dfr\) is an abelian ideal, and \(\dfr \cong \Dfr/\dfr^*\) is a Lie algebra with bracket
    \begin{align*}
        [X_1 + Y_1, X_2 + Y_2]_{\dfr} = [X_1, X_2]_{\gfr} + \alpha_{X_1} Y_2 - \beta_{Y_2} X_1 + \beta_{Y_1} X_2 - \alpha_{X_2} Y_1 + [Y_1, Y_2]_{\hfr}.
    \end{align*}
    The coadjoint action of \(\dfr\) is given by
    \begin{align}\label{eq: coadjoint action of double lie alg}
        \begin{split}
        \la \ad^{\flat}_{X_0 + Y_0} (\eta + \xi), X_1 + Y_1 \ra
        &= - \la \eta + \xi, [X_0, X_1]_{\gfr} + \alpha_{X_0} Y_1 - \beta_{Y_1} X_0 + \beta_{Y_0} X_1 - \alpha_{X_1} Y_0 + [Y_0, Y_1]_{\hfr} \ra \\
        &= \la \ad^{\flat}_{X_0} \xi + \beta_{Y_0}^{\flat} \xi - (\alpha^{\flat}_{-} \eta)^*(Y_0), X_1\ra + \la \ad^{\flat}_{Y_0} \eta + \alpha^{\flat}_{X_0} \eta - (\beta^{\flat}_{-} \xi)^*(X_0), Y_1 \ra.
        \end{split}
    \end{align}
    Comparing with \zcref{eq: bracket on Dfr} makes it clear that the action of \(\dfr\) on \(\dfr^*\) is the coadjoint action.
    The \(2\)-cocycle \(\omega\) which determines the extension \(\dfr^*_{\ab} \to \Dfr \to \dfr\) can now be computed explicitly:
    \begin{align*}
        \omega(X_1 + Y_1, X_2 + Y_2) = \vphi(X_1, X_2) + \psi(Y_1, Y_2) &- \vphi(X_1, -)^*(Y_2) + \psi(Y_2, -)^*(X_1) \\
        &+ \vphi(X_2, -)^*(Y_1) - \psi(Y_1, -)^*(X_2).
    \end{align*}

    For the converse, we define \(\gfr = \Gfr / (\dfr^* \cap \Gfr) \subseteq \dfr\) and \(\hfr = \Hfr / (\dfr^* \cap \Hfr) \subseteq \dfr\).
    Since \(\Gfr\) and \(\Hfr\) are isotropic Lie subalgebras, we must have \(\Gfr = \gfr \ltimes_{\alpha^{\flat},\omega|_{\gfr}}\hfr^*\) and \(\Hfr = \hfr \ltimes_{\beta^{\flat},\omega|_{\hfr}} \gfr^*\) for some Lie algebra representations \(\alpha : \gfr \to \gl(\hfr)\) and \(\beta : \hfr \to \gl(\gfr)\).
    Using \zcref{eq: coadjoint action of double lie alg}, we can identify \(\alpha\) and \(\beta\) with the infinitesimal dressing transformations.
    Therefore \(\Gfr\) and \(\Hfr\) are Lie bialgebras of translation type which are dual to each other.
\end{proof}

\begin{theorem}\label{th: double lie algs and lie bialgs of translation type}
    Consider a double Lie algebra \(\dfr = \gfr + \hfr\) together with matched cocycles \((\vphi, \psi)\), and write \(\alpha, \beta\) for the infinitesimal dressing transformations.
    Then \(\Gfr = \gfr \ltimes_{\alpha^{\flat},\vphi} \hfr^*_{\ab}\) is a Lie bialgebra of translation type with dual \(\Hfr = \hfr \ltimes_{\beta^{\flat},\psi} \gfr^*_{\ab}\).
    We call \(\Gfr\) and \(\Hfr\) the \emph{semi-classical limits} of the double Lie algebra \(\dfr = \gfr + \hfr\) with cocycles \((\vphi, \psi)\).

    Conversely if \(\Gfr = \gfr \ltimes_{\alpha^{\flat},\vphi} \hfr^*_{\ab}\), \(\Hfr = \hfr \ltimes_{\beta^{\flat},\psi} \gfr^*_{\ab}\) are Lie bialgebras of translation type dual to each other, then the bracket 
    \begin{align}\label{eq: bracket on dfr}
        [X_1 + Y_1, X_2 + Y_2]_{\dfr} = [X_1, X_2]_{\gfr} + \alpha_{X_1} Y_2 - \beta_{Y_2} X_1 + \beta_{Y_1} X_2 - \alpha_{X_2} Y_1 + [Y_1, Y_2]_{\hfr},
    \end{align}
    turns \(\dfr = \gfr + \hfr\) into a double Lie algebra with matched cocycles \((\vphi, \psi)\).
    These two constructions are inverses to each other.
\end{theorem}
\begin{proof}
    Let \(\dfr = \gfr + \hfr\) be a double Lie algebra with matched cocycle \((\vphi, \psi)\).
    As we have said, the data of matched cocycles \((\vphi, \psi)\) is equivalent to the data of a \(3\)-cocycle \(\Theta \in Z^3(\dfr, \R)\) whose restrictions \(\Theta|_\gfr\) and \(\Theta|_\hfr\) vanish.
    The \(3\)-cocycle \(\Theta\) is completely determined by \zcref{eq: 3 cocycle from matched cocycles}.
    We will construct a Manin triple \((\Dfr = \dfr \ltimes_{\ad^{\flat},\omega} \gfr^*_{\ab}, \Gfr, \Hfr)\) and apply \zcref{prop: double of lie bialg of translation type} to get Lie bialgebras of translation type.
    We define a \(2\)-cochain \(\omega : \dfr \wedge \dfr \to \dfr^*\) by \(\Theta(D_1, D_2, D_3) = \la \omega(D_1, D_2), D_3\ra\).
    By computing the Chevalley--Eilenberg coboundary of \(\Theta\), we obtain
    \begin{align}\label{eq: 2 and 3 cocycles}
        \begin{split}
        \delta \Theta(D_0, D_1, D_2, D_3) &= \sum_{\substack{
                \{0,1,2\} = \{i,j,k\} \\
                i < j
            }} (-1)^{i + j} \Theta([D_i, D_j], D_k, D_3)
            +  (-1)^{k + 1} \Theta(D_i, D_j, [D_k, D_3]) \\
            &= \sum_{\substack{
                \{0,1,2\} = \{i,j,k\} \\
                i < j
            }} (-1)^{i + j} \la \omega([D_i, D_j], D_k), D_3 \ra
            + (-1)^{k} \la \ad^{\flat}_{D_k} \omega(D_i, D_j), D_3 \ra \\
            &= \la \delta \omega(D_0, D_1, D_2), D_3 \ra.
        \end{split}
    \end{align}
    Therefore \(\omega\) is a \(2\)-cocycle.
    Therefore we can define the extension \(\Dfr = \dfr \ltimes_{\ad^{\flat},\omega} \dfr^*_{\ab}\).
    We note that the inner product \(\la-, -\ra_{\Dfr}\) is \(\Dfr\)-invariant:
    \begin{align*}
        &\la [D_0 + \phi_0, D_1 + \phi_1]_{\Dfr}, D_2 + \phi_2\ra_{\Dfr}
        = \la [D_0, D_1]_{\dfr}, \phi_2\ra
                + \la \ad_{D_0}^{\flat} \phi_1
                - \ad^{\flat}_{D_1} \phi_0 + \omega(D_0, D_1), D_2 \ra \\
        &\qquad= -\la D_1, \ad^{\flat}_{D_0} \phi_2 \ra
                - \la \phi_1, [D_0, D_2]_{\dfr} \ra
                + \la \ad_{D_2}^{\flat} \phi_0, D_1 \ra 
                - \Theta(D_0, D_2, D_1)\\
        &\qquad= - \la D_1 + \phi_1, [D_0 + \phi_0, D_2 + \phi_2]_{\Dfr}\ra_{\Dfr}.
    \end{align*}
    Hence it only remains to show that \(\Gfr = \gfr + \hfr^*\) and \(\Hfr = \hfr + \gfr^*\) are Lie subalgebras of \(\Dfr\).
    We only write the proof for \(\Gfr\), and the result follows by symmetry.
    Since \(\Theta|_{\gfr} = 0\), we see that \(\omega|_{\gfr}\) has image in \(\hfr^*\).
    In particular, \([X_1, X_2] \in \gfr + \hfr^*\) for all \(X_1, X_2 \in \gfr\).
    Moreover \(\hfr^*\) is a Lie subalgebra of \(\Dfr\), and it only remains to prove that \([X, \eta]_{\Dfr} \in \hfr^*\) for all \(X \in \gfr, \eta \in \hfr^*\).
    We compute
    \begin{align*}
        \la [X, \eta]_{\Dfr}, X' + Y \ra
        = - \la \eta, [X, X' + Y]_{\dfr} \ra
        = - \la \eta, \pi_{\hfr}([X, X' + Y])_{\dfr} \ra
        = - \la \eta, \pi_{\hfr}([X, Y]_{\dfr}) \ra
    \end{align*}
    for all $X, X' \in \gfr, Y \in \hfr, \eta \in \hfr^*$.
    Thus $[X, \eta]_{\Dfr} \in \hfr^*$.

    For the converse, suppose we have Lie bialgebras of translation type \(\Gfr = \gfr \ltimes_{\alpha^{\flat},\vphi} \hfr^*_{\ab}\) and \(\Hfr = \hfr \ltimes_{\beta^{\flat},\psi} \gfr^*_{\ab}\) which are dual to each other.
    We know from \zcref{prop: double of lie bialg of translation type} that their double yields the Manin triple \((\Dfr = \dfr \ltimes_{\ad^*,\omega} \dfr^*, \Gfr, \Hfr)\), where \(\dfr = \gfr + \hfr\) is a Lie algebra with bracket determined by \zcref{eq: bracket on dfr}.
    By the \(\Dfr\)-invariance of the inner product, we see that
    \begin{align*}
        \la \omega(D_0, D_1), D_2 \ra + \la D_1, \omega(D_0, D_2) \ra = 0.
    \end{align*}
    The map \(\Theta : \dfr \times \dfr \times \dfr \to \R : \Theta(D_0, D_1, D_2) := \la \omega(D_0, D_1), D_2\ra\)
    is then completely anti-symmetric.
    Moreover it is a \(3\)-cocycle thanks to \zcref{eq: 2 and 3 cocycles}.
    Finally observe that \(\omega|_{\gfr}\) has image in \(\hfr^*\) and \(\omega|_{\hfr}\) has image in \(\gfr^*\) since \(\gfr + \hfr^*\) and \(\hfr + \gfr^*\) are Lie subalgebras of \(\Dfr\).
    In particular, \(\Theta|_{\gfr} = 0\) and \(\Theta|_{\hfr} = 0\).
\end{proof}

\begin{remark}\label{rk: matched pairs and Poisson structures}
    An immediate consequence of the above theorem is that a matched pair of Lie groups \(D = G \bowtie H\) gives rise to Poisson--Lie structures on the semi-direct products \(\Gcal = G \ltimes \hfr^*\) and \(\Hcal = H \ltimes \gfr^*\) if some integrability conditions are satisfied.
    We expect smooth matched cocycles \((\Ucal, \Vcal)\) to provide a Poisson--Lie structure on some extension \(\hfr^* \to \Gcal \to G\).
    However this requires us to properly make sense of the van Est map for double Lie groups, and we will not dwell on these questions here.
    Nevertheless this suggests we should interpret the cocycle bicrossed product \(\hat G \bicrossl_{\alpha,\Ucal}^{\beta,\Vcal} H\) as a quantization of an extension \(\hfr^* \to \Gcal \to G\), which will be the focus of the next section.
\end{remark}

We conclude this discussion by describing a class of coboundary structures give rise to Lie bialgebras of translation type.
\begin{proposition}\label{prop: coboundary structures of translation type}
    Consider a Lie algebra extension \(\Gfr = \gfr \ltimes_{\lambda,\vphi} V_{\ab}\) together with \(r \in \Gfr \tensor \Gfr\) such that \((\Gfr, \delta r)\) is a Lie bialgebra.
    If \(r\) decomposes as
    \begin{align}\label{eq: coboundary of translation type}
        r = r_{\gfr,V} + r_{V,\gfr} + r_{V,V}, &&
        r_{\gfr,V} \in \gfr \tensor V,\  r_{V,\gfr} \in V \tensor \gfr,\  r_{V,V} \in V \tensor V,
    \end{align}
    then \((\Gfr, \delta r)\) is a Lie bialgebra of translation type.
\end{proposition}
\begin{proof}
    We need to show that the dual \(\Gfr^*\) decomposes as an extension \(\gfr^*_{\ab} \to \Gfr^* \to V^*\).
    Set \(r_+ := r\) and \(r_- := -r_{21}\) as usual, so that
    \begin{align*}
        [\phi_1, \phi_2]_{\Gfr^*} = \ad^{\flat}_{r_+^{\sharp} \phi_1} \phi_2 - \ad^{\flat}_{r_-^{\sharp} \phi_2} \phi_1 && \text{for all } \phi_1, \phi_2 \in \Gfr^*,
    \end{align*}
    see \cite[equation (8.8)]{majidFoundationsQuantumGroup1995}.
    We have already computed the coadjoint action of extensions in \zcref{eq: coadjoint action of extension}:
    \begin{align*}
        \ad^{\flat}_{X + v}(\xi + f) = \ad^{\flat}_X \xi - \vphi(X, -)^*(f) + (\lambda_{-} v)^* f + \lambda^{\flat}_X f
    \end{align*}
    for all \(X \in \gfr, v \in V, \xi \in \gfr^*, f \in V^*\).
    The first three terms lie in \(\gfr^*\) while the last term lies in \(V^*\).
    From the decomposition \zcref[noname]{eq: coboundary of translation type}, we observe that \(r_{\pm}^{\sharp}(\gfr^*) \subseteq V\).
    In particular we immediately have \([\gfr^*, \gfr^*]_{\Gfr^*} = 0\).
    For \(f \in V^*\), \(\xi \in \gfr^*\), we get
    \begin{align*}
        [f, \xi]_{\Gfr} = \ad^{\flat}_{r_{V,\gfr}^{\sharp} f} \xi + \ad^{\flat}_{r_{V,V}^{\sharp} f} \xi - \ad^{\flat}_{r^{\sharp}_- \xi} f
        = \ad^{\flat}_{r_{V,\gfr}^{\sharp} f} \xi - (\lambda_{-}(r_{-}^{\sharp} \xi))^*(f) \in \gfr^*.
    \end{align*}
    Hence \(\gfr^*\) is an abelian ideal of \(\Gfr^*\), which concludes.
    For the sake of completeness we compute the bracket of \(f_1, f_2 \in V^*\):
    \begin{align*}
        [f_1, f_2]_{\Gfr}
        &= \ad^{\flat}_{r_{V,\gfr}^{\sharp} f_1} f_2 + \ad^{\flat}_{r_{V,V}^{\sharp} f_1} f_2
         - \ad^{\flat}_{(r_{\gfr,V})_{21}^{\sharp} f_2} f_1 + \ad^{\flat}_{(r_{V,V})_{21}^{\sharp} f_2} f_1 \\
        &= \lambda^{\flat}_{r_{V,\gfr}^{\sharp} f_1} f_2 + \lambda^{\flat}_{(r_{\gfr,V})_{21}^{\sharp} f_2} f_1
         - \vphi\left(r_{V,\gfr}^{\sharp} f_1, -\right)^*(f_2)
         - \vphi\left((r_{\gfr,V})_{21}^{\sharp} f_2, -\right)^*(f_1)\\
        &\qquad\qquad+ \left(\lambda_{-} r_{V,V}^{\sharp} f_1\right)^*(f_2) + \left(\lambda_{-} (r_{V,V})_{21}^{\sharp} f_2\right)^*(f_1).
    \end{align*}
\end{proof}

\section{Abelian approximations and the Poisson--Fourier transform}\label{section: abelian approximation}
\zcref[S]{th: double lie algs and lie bialgs of translation type} suggest that a cocycle bicrossed product \(\hat G \bicrossl_{\alpha,\Ucal}^{\beta,\Vcal} H\) should be interpreted as a quantization of an extension \(V \to \Gcal \to G\).
To do so one we want to compare the von Neumann algebras \(W^*(G) \ltimes_{\alpha,\Ucal} \Lrm^{\infty}(H)\) with \(W^*(\Gcal)\).
Note that both \(\Lrm^{\infty}(H) \leq W^*(G) \ltimes_{\alpha,\Ucal} \Lrm^{\infty}(H)\) and \(W^*(V) \leq W^*(\Gcal)\) are subalgebras.
To compare the two we may make use of the Fourier transform \(\Fcal_V : W^*(V) \cong \Lrm^{\infty}(\hat V)\) as well as a \(G\)-equivariant measure class isomorphism \(\eta_H : H \to \hat V\) which is well-behaved with respect to \(\Ucal\), see \zcref{eq:abApproxCocycle}.
The composition \(\Fcal_H := \Fcal_{\hat V} \eta_{H*} : W^*(G) \ltimes_{\alpha,\Ucal} \Lrm^{\infty}(H) \to W^*(\Gcal)\) is then an isomorphism of von Neumann algebras.
Such an \(\eta_H\) will be called an \emph{abelian approximation of \(H\)}, and \(\Fcal_H\) the \emph{Poisson--Fourier \(\frac{1}{2}\)-transform} associated with \(\eta_H\).
We can now quantize \((W^*(\Gcal), \hat \Delta_\Gcal)\) by replacing the coproduct \(\hat \Delta_{\Gcal}\) with \(\hat \Delta^{\beta,\Vcal}_{\Gcal} := (\Fcal_H \tensor \Fcal_H) \Delta^{\beta,\Vcal} \Fcal_H^{-1}\).
Suppose that \(\eta_G : G \to \hat W\) is an abelian approximation of \(G\) for a suitable extension \(W \to \Hcal \to H\).
Then by making use of \zcref{prop:dualityBicrossed} we obtain an isomorphism of \lcqgs
\begin{align*}
    \Fcal_{\Gcal} := \Fcal_G \Sigma \Fcal_H^{-1} : \left(
                W^*(\Gcal),\hat \Delta_{\Gcal}^{\beta,\Vcal}
            \right)\sphat\,
            \xrightarrow{\cong}
            \left(
                W^*(\Hcal), \hat \Delta_{\Hcal}^{\alpha,\Ucal}
            \right).
\end{align*}
We call \(\Fcal_{\Gcal}\) the \emph{Poisson--Fourier transform} between \(\Gcal\) and \(\Hcal\) and interpret the above equation as implement the \qdp, see \zcref{th:PFT
}
\begin{definition}\label{def:abelianApprox}
    Let $G \bowtie H = D$ be a matched pair with dressing transformations $(\alpha, \beta)$ and cocycles $(\Ucal, \Vcal)$.
    A \emph{measurable abelian approximation of $H$} consists of a \lcab $V$ equipped with a right $G$-action $\rho : V \times G \to V$ by continuous automorphisms and a Borel $2$-cocycle $\psi : G \times G \to V$, together with a measure class isomorphism $\eta_H : H \to \hat V$ such that for almost every $g, g' \in G, h \in H, v \in V$, we have
    \begin{align}\label{eq:abApproxAction}
        \ppairing{\rho(v)_g}{\eta_H(h)} = \ppairing{v}{\eta_H(\alpha_g(h))},
    \end{align}
    and
    \begin{align}\label{eq:abApproxCocycle}
        \ppairing{\psi(g,g')}{\eta_H(h)} = \Ucal(g, g', h).
    \end{align}
    We refer to the extension $\Gcal = G \ltimes_{\rho,\psi} V$ as the \emph{classical limit of $\hat G \bicrossl_{\alpha,\Ucal}^{\beta,\Vcal} H$} associated with the abelian approximation $\eta_H : H \to \hat V$.

    By a \emph{measurable abelian approximation of $G$} for the matched pair $G \bowtie H = D$ with cocycles $(\Ucal,\Vcal)$, we mean a measurable abelian approximation of $G$ for the dual matched pair $H \bowtie G = D$ with cocycles $(\Vcal,\Ucal)$.
\end{definition}

\begin{remark}
    Any measure class isomorphism $\eta_0 : H \to \hat V$ uniquely determine a Borel action $\rho_0 : V \times G \to V$ and a Borel $2$-cocycle $\psi_0 : G \times G \to V$ using \zcref{eq:abApproxAction,eq:abApproxCocycle}.
    They key point of the above definition is that \(\rho\) is not just a Borel action, but a continuous action by automorphisms.
\end{remark}

Clearly, an abelian approximation $\eta_H : H \to \hat V$ induces an isomorphism of crossed product von Neumann algebras
\[
    {\eta_H}_* : W^*(G) \ltimes_{\alpha,\Ucal} \Lrm^{\infty}(H) \to W^*(G) \ltimes_{\hat \rho, \psi} \Lrm^{\infty}(\hat V),
\]
where $\hat \rho : G \times \hat V \to \hat V$ is the left $G$-action dual to $\rho$.

\begin{construction}[The Poisson--Fourier $\frac{1}{2}$-transform]\label{constr:PFHalfTransform}
    Consider a matched pair $G \bowtie H = D$ with dressing transformations $(\alpha, \beta)$ and cocycles $(\Ucal, \Vcal)$, as well as a measurable abelian approximation $\eta_H : H \to \hat V$ with classical limit $\Gcal = G \ltimes_{\rho,\psi} V$.
    We note that \(\id \tensor \Fcal_V^* : W^*(G) \ltimes_{\hat \rho,\psi} \Lrm^{\infty}(\hat V) \to W^*(\Gcal)\) is an isomorphism.
    In particular, the unitary
    \[
        \Fcal_H = \Fcal_V^* {\eta_H}_* : \Lrm^2(G) \tensor \Lrm^2(H) \to \Lrm^2(\Gcal) = \Lrm^2(G) \tensor \Lrm^2(V)
    \]
    provides an isomorphism of von Neumann algebras
    \[
        \Fcal_H : W^*(G) \ltimes_{\alpha,\Ucal} \Lrm^{\infty}(H) \cong W^*(\Gcal).
    \]
    We transport the coproduct $\Delta^{\beta,\Vcal}$ and the Haar weights on $\hat G \bicrossl_{\alpha,\Ucal}^{\beta,\Vcal} = (W^*(G) \ltimes_{\alpha,\Ucal} \Lrm^{\infty}(H), \Delta^{\beta,\Vcal})$ to a new coproduct $\hat \Delta_{\Gcal}^{\beta,\Vcal}$ and new Haar weights on $W^*(\Gcal)$, so that we get an isomorphism of \lcqgs
    \[
        \Fcal_H : \hat G \bicrossl_{\alpha,\Ucal}^{\beta,\Vcal} H \cong (W^*(\Gcal), \hat \Delta^{\beta,\Vcal}_{\Gcal}).
    \]
    We interpret $(W^*(\Gcal), \hat \Delta^{\beta,\Vcal}_{\Gcal})$ as a quantization of $(W^*(\Gcal), \hat \Delta_{\Gcal})$, and we call $\Fcal_H$ the \emph{Poisson--Fourier $\frac{1}{2}$-transform} associated with the abelian approximation $\eta_H : H \to \hat V$.
\end{construction}

\begin{remark}\label{rk: semi-classical limit}
    One can ask how legitimate it is to interpret $(W^*(\Gcal),\hat \Delta^{\beta,\Vcal}_{\Gcal})$ as a quantization of $(W^*(\Gcal),\hat \Delta_{\Gcal})$.
    Because of the generality we are working in, it is not straightforward to include a deformation parameter in the construction.
    Moreover, any attempt at performing the semi-classical analysis of this sort of quantization would require additional regularity conditions on the abelian approximation \(\eta_H\).
    We still believe this is reasonable for the following reasons:
    If we are working with Lie groups over a local field $k$ and we take $V$ to be a $k$-vector space, then from $\eta_H : H \to \hat V$ we can define
    \[
        \eta_{H,\hbar} : H \to \hat V : \eta_{H,\hbar}(h) = \frac{1}{\hbar} \eta_H(h)
    \]
    for $\hbar \in k \setminus \{0\}$.
    In many examples this yields a continuous family $(W^*(\Gcal), \hat \Delta^{\hbar}_{\Gcal})$ which converges to $(W^*(\Gcal), \hat \Delta_{\Gcal})$ as $\hbar \to 0$, and the resulting products \(\star_{\hbar}\) on the Fourier algebra \(\Acal(\Gcal)\) quantize the expected Poisson structure.

    Moreover we will show in an upcoming paper that the bicrossed product of the double Lie group associated to a triangular (resp. quasi-triangular) Lie bialgebra of translation type through \zcref{th: double lie algs and lie bialgs of translation type} is itself a triangular (resp. quasi-triangular) \lcqg.
    Important properties of the Poisson--Lie structure of \(\Gcal\) are therefore reflected in the quantization \((W^*(\Gcal),\hat \Delta^{\beta,\Vcal}_{\Gcal})\).
\end{remark}

We now turn our attention to the quantum duality principle.
Consider two extensions $\Gcal = G \ltimes_{\rho,\psi} V$ and $\Hcal = H \ltimes_{\sigma,\omega} W$ of \lcgs, where $V, W$ are abelian, $\rho : V \times G \to V$, $\sigma : W \times H \to W$ are right actions by continuous automorphisms, and $\psi : G \times G \to V$, $\omega : H \times H \to W$ are Borel $2$-cocycles.
In practice $\Gcal$ and $\Hcal$ will often be dual Poisson--Lie groups.
We want to make use of the theory of abelian approximations developed in the previous section to quantize both $\Gcal$ and $\Hcal$, and combine the resulting Poisson--Fourier $\frac{1}{2}$-transforms to realise the quantum duality principle.
For this purpose, suppose we have measure class isomorphisms $\eta_H : H \to \hat V$ and $\eta_G : G \to \hat W$.
We define Borel actions $\alpha : G \times H \to H$ and $\beta : H \times G \to G$ by setting for almost every $g \in G, h \in H, v \in V, w \in W$
\begin{equation}\label{eq:matchedActionsFromAbApprox}
\begin{split}
    \ppairing{v}{\eta_H(\alpha_g(h))} &= \ppairing{\rho(v)_g}{\eta_H(h)}; \\
    \ppairing{w}{\eta_G(\beta_h(g))} &= \ppairing{\sigma(w)_h}{\eta_G(g)}.
\end{split}
\end{equation}
Similarly, we obtain Borel maps $\Ucal : G \times G \times H \to U(1)$ and $\Vcal : H \times H \times G \to U(1)$ using the equations
\begin{equation}\label{eq:matchedCocyclesFromAbApprox}
\begin{split}
    \Ucal(g_1, g_2, h) &= \ppairing{\psi(g_1, g_2)}{\eta_H(h)}; \\
    \Vcal(h_1, h_2, g) &= \ppairing{\omega(h_1,h_2)}{\eta_G(g)}
\end{split}
\end{equation}
for almost every $g, g_1, g_2 \in G$, $h, h_1, h_2 \in H$.
One should compare the above equations with \zcref{eq:abApproxAction,eq:abApproxCocycle}.

\begin{definition}\label{def:matchedAbApprox}
    In the above setting, we say that the Borel maps $(\eta_H, \eta_G)$ form \emph{matched abelian approximations} if they satisfy the following axioms:
    \begin{enumerate}
        \item The map $\tau : G \times H \to G \times H : (g, h) \mapsto (\beta_h(g), \alpha_g(h))$ is a measure class isomorphism.
        \item For almost every $g, g' \in G$, $h, h' \in h$, we have
        \begin{align}\label{eq:matchedAbApproxActions}
            \alpha_g(h' h) = \alpha_{\beta_h(g)}(h') \alpha_g(h);
            && \beta_h(g' g) = \beta_{\alpha_g(h)}(g') \beta_h(g).
        \end{align}
        \item For almost every $g_1, g_2 \in G$, $h_1, h_2 \in H$, we have
        \begin{equation}\label{eq:matchedAbApproxCocycles}
        \begin{split}
            \Vcal(h_1, h_2, g_1 g_2) \cj{\Ucal}(g_1, g_2, h_1 h_2)
            ={} & \Vcal(\alpha_{\beta_{h_2}(g_2)}(h_1, \alpha_{g_2}(h_2), g_1))
                \Vcal(h_1, h_2, g_2) \\
            &{}\cj{\Ucal}(\beta_{\alpha_{g_2}(h_2)}(g_1), \beta_{h_2}(g_2), h_1)
                \cj{\Ucal}(g_1, g_2, h_2).
        \end{split}
        \end{equation}
    \end{enumerate}
\end{definition}

\begin{remark}
    The \zcref{eq:matchedAbApproxActions,eq:matchedAbApproxCocycles} come from the identities related to matched pairs and their cocycles, see \zcref{eq:dresTransIds,eq:matchedCocycles}.
\end{remark}

As expected, matched abelian approximations are measurable abelian approximations in the sense of \zcref{def:abelianApprox}.
\begin{theorem}\label{th:matchAbApproxAreAbApprox}
    Consider extensions $\Gcal = G \ltimes_{\rho,\psi} V$ and $\Hcal = H \ltimes_{\sigma,\omega} W$ as above, and matched abelian approximations $\eta_H : H \to \hat V$ and $\eta_G : G \to \hat W$.
    We define $(\alpha,\beta)$ and $(\Ucal,\Vcal)$ using \zcref{eq:matchedActionsFromAbApprox,eq:matchedCocyclesFromAbApprox}.

    Then there exists a \lcg $D$ which admits $G, H$ as closed subgroups such that $G \bowtie H = D$ is a matched pair with dressing transformations $(\alpha,\beta)$ and cocycles $(\Ucal,\Vcal)$.
    The Borel maps $\eta_H$ and $\eta_G$ are measurable abelian approximations of $H$ and $G$ respectively for the matched pair $G \bowtie H = D$ with cocycles $(\Ucal,\Vcal)$, in the sense of \zcref{def:abelianApprox}.
    The classical limit associated to $\eta_H$ is $\Gcal$, and the classical limit associated to $\eta_G$ is $\Hcal$.
\end{theorem}
\begin{proof}
    From eqs. \ref{eq:matchedAbApproxActions} we see that the $*$-isomorphism
    $\tau^* : \Lrm^{\infty}(G) \htensor \Lrm^{\infty}(H) \to \Lrm^{\infty}(G) \htensor \Lrm^{\infty}(H)$
    determined by
    $\tau^*(\vphi)(g,h) = \vphi(\beta_h(g), \alpha_g(h))$
    is a matching of $\Lrm^{\infty}(G)$ and $\Lrm^{\infty}(H)$ with trivial cocycles in the sense of \cite[Definition 2.1]{vaesExtensionsLocallyCompact2003}.
    Hence we can perform the bicrossed product construction and obtain a \lcqg $\hat G \bicrossl_{\alpha}^{\beta} H$.
    Its multiplicative unitary comes from a pentagonal transformation $v : X \times X \to X \times X$ on $X = G \times H$, see \cite[Section 5.1]{baajNonSemiRegularQuantumGroups2003}.
    It is given by the formula
    \[
        v(g_1, h_1, g_2, h_2) = (\beta_{\alpha_{g_1}(h_1)^{-1} h_2}(g_2) g_1, h_1, g_2, \alpha_{g_1}(h_1)^{-1} h_2).
    \]
    By \cite[Proposition 5.1]{baajNonSemiRegularQuantumGroups2003} we conclude that this pentagonal transformation $v$ is produced from a matched pair $G \bowtie H = D$ of \lcgs.
    We note that we use slightly different conventions compared to \cite{baajNonSemiRegularQuantumGroups2003}.
    We can easily translate our setting in theirs using the group inverse on $H$, see the discussion under \cite[Definition 3.3]{baajNonSemiRegularQuantumGroups2003}.

    The \zcref{eq:matchedAbApproxCocycles} is precisely the necessary condition for $(\Ucal,\Vcal)$ to form matched cocycles for the matched pair $G \bowtie H = D$.
    Finally, \zcref{eq:matchedActionsFromAbApprox,eq:matchedCocyclesFromAbApprox} correspond exactly to the \zcref{eq:abApproxAction,eq:abApproxCocycle} which characterise abelian approximations.
\end{proof}
\begin{remark}
    It is worth noting that the ambient group $D$ in the above theorem can also be described as the double crossed product of $G$ and $H$ for the matching $\tau^*$ \cite{baajDoubleCrossedProducts2005}.
\end{remark}

We combine \zcref{th:matchAbApproxAreAbApprox} with \zcref{constr:PFHalfTransform} to obtain quantizations $(W^*(\Gcal), \hat\Delta^{\beta,\Vcal}_{\Gcal})$ and $(W^*(\Hcal),\hat \Delta^{\alpha,\Ucal}_{\Hcal})$ of $\Gcal$ and $\Hcal$ respectively, together with the corresponding Poisson--Fourier $\frac{1}{2}$-transforms
\begin{align*}
    \Fcal_H : \hat G \bicrossl_{\alpha,\Ucal}^{\beta,\Vcal} H \xrightarrow{\cong} \left(W^*(\Gcal),\hat\Delta_{\Gcal}^{\beta,\Vcal}\right); &&
    \Fcal_G : \hat H \bicrossl_{\beta,\Vcal}^{\alpha,\Ucal} G \xrightarrow{\cong} \left(W^*(\Hcal),\hat\Delta_{\Hcal}^{\alpha,\Ucal}\right).
\end{align*}
By \zcref{prop:dualityBicrossed}, the flip  map $\Sigma : \Lrm^2(G) \tensor \Lrm^2(H) \to \Lrm^2(H) \tensor \Lrm^2(G)$ gives an isomorphism of \lcqgs
\[
    \Sigma : \left(\hat G \bicrossl_{\alpha,\Ucal}^{\beta,\Vcal} H\right)\sphat\, \xrightarrow{\cong} \hat H \bicrossl_{\beta,\Vcal}^{\alpha,\Ucal} G.
\]
By putting all of this together we realise the quantum duality principle, as summarised in the following theorem:

\begin{theorem}[Quantum duality principle]\label{th:PFT}
    In the above setting, consider the unitary operator
    \begin{equation}
        \Fcal_{\Gcal} = \Fcal_G \Sigma \Fcal_H^* : \Lrm^2(\Gcal) \to \Lrm^2(\Hcal).
    \end{equation}
    Then $\Fcal_{\Gcal}$ provides an isomorphism of \lcqgs
    \begin{equation}
        \Fcal_{\Gcal} : \left(
                W^*(\Gcal),\hat \Delta_{\Gcal}^{\beta,\Vcal}
            \right)\sphat\,
            \xrightarrow{\cong}
            \left(
                W^*(\Hcal), \hat \Delta_{\Hcal}^{\alpha,\Ucal}
            \right).
    \end{equation}
    In this sense $\Fcal_{\Gcal}$ implements the quantum duality principle, and we call it the \emph{Poisson--Fourier transform} between $\Gcal$ and $\Hcal$.
\end{theorem}

\section{Abelian approximations from coboundary structures}\label{section: abelian approximation from coboundary structure}
We now describe a general strategy to construct abelian approximations from Lie bialgebras of translation type with coboundary structures.
This works especially well when the \(r\)-matrix is triangular and has some applications when it is quasi-triangular, as explained in \zcref{cor: abelian approximation for triangular TT,cor: abelian approximation for quasi-triangular TT}.
For simplicity we restrict ourselves to the case without cocycles, but we hope some of the issues with cocycles can be resolved using the methods of \cite{kalyuzhnyiFindingCocyclesBicrossed2007}.
Consider a double Lie algebra \(\dfr = \gfr + \hfr\) and write \(\alpha\), \(\beta\) for the infinitesimal dressing transformations.
From \zcref{th: double lie algs and lie bialgs of translation type} we know that \(\Gfr = \gfr \ltimes_{\alpha^{\flat}} \hfr^*_{\ab}\) is a Lie bialgebra of translation type with dual \(\Hfr = \hfr \ltimes_{\beta^{\flat}} \gfr^*_{\ab}\).
Furthermore we assume this double Lie algebra integrates to a matched pair of Lie groups \(D = G \bowtie H\).
We wish to construct an abelian approximation \(\eta_H : H \to \hfr\).
In many situations we can do so many making use of a \emph{reductive complement} for \(\gfr\).

\begin{definition}\label{def: reductive complement}
    Let \(\dfr = \gfr + \hfr\) be a double Lie algebra.
    A \emph{reductive complement for \(\gfr\)} consists of a linear subspace \(\pfr \subseteq \dfr\) such that \(\dfr = \gfr \oplus \pfr\) as a vector space, and
    \begin{align*}
        [\gfr, \pfr] \subseteq \pfr.
    \end{align*}

    A double Lie algebra with reductive complement \(\dfr = \gfr + \hfr = \gfr + \pfr\) is said to be \emph{integrable} if there exists a matched pair of Lie groups \(D = G \bowtie H\) integrating the Lie algebras \(\dfr = \gfr + \hfr\) such that
    \begin{enumerate}
        \item\label{axiom: G-equivariance} $\Ad_g(\pfr) \subseteq \pfr$ for all $g \in G$;
        \item\label{axiom: exponential} The restricted exponential map $\exp_{\pfr} := \exp_{D}|_{\pfr}: \pfr \to P := \exp_D(\pfr)$ is a diffeomorphism onto its image;
        \item\label{axiom: total decomposition} We have \(D = GP\) and \(G \cap P = \{e_D\}\).
    \end{enumerate}
    In this case we say that \(D = G \bowtie H = G P\) is a double Lie group with exponential reductive complement.
\end{definition}

\begin{example}
    Let \(D\) be a connected semi-simple Lie group.
    Then the Iwasawa and Cartan decompositions \(D = K \bowtie AN = K P\) form a double Lie group with exponential reductive complement. 
\end{example}

Consider a double Lie group with exponential reductive complement \(D = G \bowtie H = GP\).
Clearly \(c_g(p) := g p g^{-1} \in P\) for all \(g \in G, p \in P\), and thus \(G\) acts on \(P\) by conjugation.
In particular, the map \(G \times P \to D : (g, p) \mapsto gp\) is a measure class isomorphism \cite[Theorem~2.58 and Theorem~2.61]{follandCourseAbstractHarmonic2016}.
Using \cite[Proposition~3.2]{baajNonSemiRegularQuantumGroups2003}, we observe that \(\chi : G \times H \to G \times P : \chi(g, h) = (g \Lrm^G_P(h), \Rrm^P(h))\) is also a measure class isomorphism.
However if \(N \subset H\) is a Borel subset of measure \(0\), then \(\chi(G \times N) = G \times \Rrm^P(N)\) has measure \(0\) and thus so does \(\Rrm^P(N) \subset P\).
Using a similar argument with the map \(\chi' : G \times P \to G \times H : \chi'(g, p) = (g \Lrm^G_H(p), \Rrm^H(p))\) we conclude that the restriction \(\Rrm^P : H \to P\) is a measure class isomorphism. 
Moreover note that for almost all \(h \in H\), \(g \in G\),
\begin{align*}
    \Rrm^P(\alpha_g(h)) = \Rrm^P(\Rrm^H(h g^{-1})) = \Rrm^P(h g^{-1}) = \Rrm^P(\Rrm^P(h) g^{-1}) = c_g(\Rrm^P(h)).
\end{align*}
We have proved the following proposition:
\begin{proposition}\label{prop: abelian approximation from lie--reductive doubles}
    Suppose that \(D = G \bowtie H = GP\) is a double Lie group with exponential reductive complement.
    Then the map
    \begin{align*}
        \eta_H = \exp_{\pfr}^{-1} \circ \Rrm^P : H \to \pfr
    \end{align*}
    is an abelian approximation for \(H\).
\end{proposition}

\begin{remark}
    We can relax the conditions \zcref[noname]{axiom: exponential} and \zcref[noname]{axiom: total decomposition} as follows:
    \begin{itemize}
        \item We only really need the decomposition \(D = GP\) to hold almost everywhere, and so we may ask instead that \(D \setminus GP\) has measure \(0\);
        \item We can replace the exponential map \(\exp|_{\pfr} : \pfr \to P\) by any measure class isomorphism \(\Exp : \pfr \to P\).
            In this case the measure on \(P\) should be taken to come from any strongly \(G\)-quasi-invariant measure on the homogeneous space \(D/G\) when it is identified with \(P\) under the projection map \(\pi : D \to D/G\).
    \end{itemize}
    It is also reasonable to ask that \(\Exp\) is smooth almost everywhere; \(\Exp(0) = e_D\); and \(\dd_0 \Exp = \id_{\pfr}\).
    
\end{remark}

At the time of writing we are not aware of a general method to determine whether a given double Lie algebra with reductive complement \(\dfr = \gfr + \hfr = \gfr + \pfr\) is integrable, as constructing the matched pair \(D = G \bowtie H\) is already non-trivial.
However is \(\pfr\) is a solvable ideal whose adjoint action has no purely imaginary eigenvalues, then its associated simply connected Lie group is exponential \cite{dixmierLapplicationExponentielleDans1957,saitoCertainsGroupesLieII1957}.
Then it only remains to find a suitable Lie subgroup \(H\) in the semi-direct product \(D = G \ltimes P\).
In many cases one can choose \(\pfr\) such that \([\pfr, \pfr]_{\dfr} \subseteq \gfr\).
One might then be able to show conditions \zcref[noname]{axiom: exponential,axiom: G-equivariance,axiom: total decomposition} using the methods of \cite{voglaireStronglyExponentialSymmetric2014}, but we will not investigate this strategy here.

We have already discussed in \zcref{prop: double of lie bialg of translation type,th: double lie algs and lie bialgs of translation type} how any Lie bialgebra of translation type \(\Gfr = \gfr \ltimes_{\alpha} \hfr^*_{\ab}\) with dual \(\Hfr = \hfr \ltimes_{\beta} \gfr^*_{\ab}\) gives rise to a double Lie algebra \(\dfr = \gfr + \hfr\) through the double construction.
Our next goal is then to construct a reductive complement \(\dfr = \gfr + \pfr\) from a coboundary structure \(r \in \gfr \tensor \gfr\) as in \zcref{prop: coboundary structures of translation type}.
We first discuss the structure of the double of a coboundary Lie bialgebra.
While the formulation of the next two theorems are novel, most if not all of the statements were already known before.
Indeed the splitting maps \(j_{\pm} : \gfr^* \to \dfr\) come from \cite[Lemma~3.1]{alekseevLinearizationPoissonLie2016}.
\zcref[S]{item: ideals from quasi-triangular r-matrix,item: abelian ideal from triangular r-matrix} of \zcref{th: dual reductive complements in double} appear in some form in \cite{kosmann-schwarzbachPoissonLieGroupsComplete1988,drinfeldAlmostCocommutativeHopf1989}, as well as \cite{weinsteinClassicalSolutionsQYBE1992} which contains the equation \([\pfr_+, \pfr_-] = 0\).
\begin{theorem}\label{th: coboundary structures and reductive complements}
    Let \((\dfr, \gfr, \gfr^*)\) be a Manin triple.
    There is a bijection between \(\gfr\)-equivariant splittings \(j : \gfr^* \to \dfr\) of the sequence \(\gfr \to \dfr \to \gfr^*\) and coboundary structures \(r \in \gfr \tensor \gfr\).
    If \(r \in \gfr \tensor \gfr\) is a coboundary structure, then the corresponding splitting is
    \begin{align*}
        &j_r : \gfr^* \to \dfr\\
        &j_r(\xi) := \xi - r^{\sharp}(\xi).
    \end{align*}
    In particular if we set \(\pfr_r := j_r(\gfr^*)\), we obtain a bijection between reductive complements \(\dfr = \gfr + \pfr\) of \(\gfr\) and coboundary structures \(r \in \gfr \tensor \gfr\).
\end{theorem}
\begin{proof}
    The first part of the statement is exactly \cite[Lemma~3.1]{alekseevLinearizationPoissonLie2016}.
    For the second part, observe that if \(\dfr = \gfr + \pfr\) is a reductive complement for \(\gfr\), then the projection \(\pi_{\pfr} : \dfr \to \pfr\) restricts to a \(\gfr\)-equivariant injective linear map \(\pi_{\pfr}|_{\gfr^*} : \gfr^* \to \pfr \subseteq \dfr\).
    One easily verifies that these two constructions are inverses to each other.
\end{proof}

\begin{theorem}\label{th: dual reductive complements in double}
    Let \((\gfr, r)\) be a coboundary Lie bialgebra and set \(r_+ := r\), \(r_- = -r_{21}\).
    We let \(\pfr_{+}, \pfr_{-}\) denote the reductive complements associated with \(r_+, r_-\) respectively.
    Then
    \begin{align*}
        [\pfr_+, \pfr_-] \subseteq \gfr.
    \end{align*}
    Moreover,
    \begin{enumerate}
        \item\label{item: symmetric space from skew-symmetric coboundary} If \(r\) is skew-symmetric, then \(\pfr_+ = \pfr_- =: \pfr\), and \([\pfr, \pfr] \subseteq \gfr\);
        \item\label{item: ideals from quasi-triangular r-matrix} If \(r\) is quasi-triangular, then \(\pfr_+, \pfr_-\) are ideals and \([\pfr_+, \pfr_-] = 0\);
        \item\label{item: abelian ideal from triangular r-matrix} If \(r\) is triangular, then \(\pfr_+ = \pfr_- =: \pfr\) is an abelian ideal and \(\dfr \cong \gfr \ltimes_{\ad^{\flat}} \gfr^*\).
    \end{enumerate}
\end{theorem}
\begin{proof}
    Define \(j_{\pm} : \gfr^* \to \dfr : \xi \mapsto \xi - r_{\pm}^{\sharp}(\xi)\) as in the previous theorem, so that \(\pfr_{\pm} = j_{\pm}(\gfr^*)\).
    For any $\xi, \zeta \in \gfr^*$, we compute
    \begin{align*}
        &[j_+(\xi), j_{-}(\zeta)]_{\dfr}
        = [\xi - r_+^{\sharp}(\xi), \zeta - r_-^{\sharp}(\zeta)]_{\dfr} \\
        &\quad= \ad^{\flat}_{r_+^{\sharp}(\xi)} \zeta - \ad^{\flat}_{r^{\sharp}_-(\zeta)} \xi - \ad^{\flat}_{\xi} r^{\sharp}_-(\zeta) + \ad^{\flat}_{r^{\sharp}_{-}(\zeta)} \xi - \ad^{\flat}_{r_+^{\sharp}(\xi)} \zeta + \ad^{\flat}_{\zeta} r_+^{\sharp}(\xi) + [r^{\sharp}_+(\xi), r^{\sharp}_-(\zeta)]_{\gfr} \\
        &\quad= [r_+^{\sharp}(\xi), r_{-}^{\sharp}(\zeta)]_{\gfr} + \ad^{\flat}_{\zeta} r^{\sharp}_+(\xi) - \ad^{\flat}_{\xi} r^{\sharp}_-(\zeta) \in \gfr.
    \end{align*}
    If \(r\) is skew-symmetric, then \(r_+ = r_-\) and \zcref[noname]{item: symmetric space from skew-symmetric coboundary} is immediate.
    Furthermore \zcref[noname]{item: abelian ideal from triangular r-matrix} follows from \zcref[noname]{item: symmetric space from skew-symmetric coboundary,item: ideals from quasi-triangular r-matrix}, and it only remains to show \zcref[noname]{item: ideals from quasi-triangular r-matrix}.

    Suppose that \(r \in \gfr \tensor \gfr\) is quasi-triangular and define \(f_{\pm} : \dfr \to \gfr\) by \(f_{\pm}(X + \xi) := X + r^{\sharp}_{\pm}(\xi)\).
    It is well-known that \(f_{\pm}\) are Lie algebra homomorphisms, see for example \cite[Proposition~2.3]{weinsteinClassicalSolutionsQYBE1992}.
    Hence \(\pfr_{\pm} = \ker f_{\pm}\) are ideals in \(\dfr\),
    and \([\pfr_+, \pfr_-]_{\dfr} \subseteq \gfr \cap \pfr_+ = 0\).
\end{proof}

Applying the above theorem to coboundary Lie bialgebras of translation type and quotienting out the abelian ideal \(\dfr^* \leq \Dfr\) gives the following corollary.
\begin{corollary}\label{cor: reductive complement from coboundary of translation type}
    Let \(\Gfr = \gfr \ltimes_{\alpha} \hfr^*_{\ab}\) be a Lie bialgebra of translation type without cocycle whose dual is \(\Hfr = \hfr \ltimes_{\beta} \gfr^*_{\ab}\).
    We write \(\dfr = \gfr + \hfr\) for the associated double Lie algebra and \(\Dfr = \dfr \ltimes_{\ad^{\flat}} \dfr^*\).
    Suppose \(\Gfr\) admits a coboundary structure \(r \in \Gfr \tensor \Gfr\) and write \(r_+ = r\), \(r_- = -r_{21}\).
    If \(\Pfr_{\pm} \subseteq \Dfr\) denote the corresponding reductive complements of \(\Gfr\) in \(\Dfr\), then \(\pfr_{\pm} := \Pfr_{\pm} / \dfr^* \subseteq \Dfr / \dfr^* \cong \dfr\) form reductive complements \(\dfr = \gfr + \pfr_{\pm}\).
    Moreover,
    \begin{enumerate}
        \item If \(r\) is skew-symmetric, then \(\pfr_+ = \pfr_- =: \pfr\) and \([\pfr, \pfr]_{\dfr} \subseteq \gfr\);
        \item If \(r\) is quasi-triangular, then \(\pfr_+, \pfr_-\) are ideals and \([\pfr_+, \pfr_-]_{\dfr} = 0\);
        \item If \(r\) is triangular, then \(\pfr_+ = \pfr_- =: \pfr\) is an abelian ideal and \(\dfr \cong \gfr \ltimes_{\alpha} \hfr_{\ab}\).
    \end{enumerate}
\end{corollary}

\begin{example}\label{example: tangent Lie bialgebra}
    Let $(\gfr, r)$ be any coboundary Lie bialgebra.
    Then $\Gfr := \gfr \ltimes_{\ad} \gfr_{\ab}$ is a coboundary Lie bialgebra of translation type with coboundary structure
    \[
        r_{\Gfr} = r_{14} + r_{23} \in \Gfr \tensor \Gfr.
    \]
    Its dual is $\Gfr^* = \gfr^* \ltimes \gfr^*_{\ab}$, and the associated double Lie algebra with reductive complement is exactly given by \zcref{th: coboundary structures and reductive complements} applied to \(\gfr\):
    \begin{align*}
        \dfr = \gfr + \gfr^* = \gfr + j_r(\gfr^*),
    \end{align*}
    where $\dfr$ denotes the double of $\gfr$.
    This suggests that the bicrossed product of dual Poisson--Lie groups \(G\), \(G^*\) should be interpreted as a quantization of the tangent Lie group \(T G\).
    We note that if $r$ is quasi-triangular or triangular, then so is $r_{\Gfr}$.
\end{example}

\begin{example}\label{example: cotangent bundle}
    Any Lie algebra $\gfr$ can be turned into a Lie bialgebra by equipping it with the trivial cocommutator $\gamma_{\gfr} = 0$.
    Its double $\Gfr = \gfr \ltimes_{\ad^{\flat}} \gfr^*_{\ab}$ is then a factorisable Lie bialgebra of translation type when equipped with the $r$-matrix
    \begin{align*}
        r = \sum_{i = 1}^n X_i \tensor X_i^*,
    \end{align*}
    where $\{X_1, \dots, X_n\} \subset \gfr$ is any basis of $\gfr$, and $\{X_1^*, \dots, X_n*\} \subset \gfr^*$ denotes the corresponding dual basis of $\gfr^*$.
    The double of a factorisable Lie bialgebra takes a very simple form \cite{semenov-tian-shanskyDressingTransformationsPoisson1985}.
    It can be described using the Manin triple
    \begin{align*}
        (\Dfr = \Gfr \times \Gfr, \diag(\Gfr), (r_+^{\sharp}, r_-^{\sharp})(\Gfr^*)),
    \end{align*}
    where $\diag(X) = (X, X)$ is the diagonal map and $(r_+^{\sharp}, r_-^{\sharp})(\xi) = (r_+^{\sharp}(\xi), r_{-}^{\sharp}(\xi))$ for $\xi \in \Gfr^*$.
    If we write \(\dfr = \gfr \times \gfr\), we obtain the double Lie algebra with reductive complement
    \begin{align*}
        \dfr = \diag(\gfr) + (0 \times \gfr).
    \end{align*}
    Here $0 \times \gfr$ plays the role of both $\hfr$ and $\pfr$ at the same time.
    Note that the decomposition \(\dfr = \diag(\gfr) + (0 \times \gfr)\) presents \(\gfr \times \gfr\) as the semi-direct product \(\dfr \cong \gfr \ltimes_{\ad} \gfr\).
\end{example}

\begin{example}
    We can slightly generalise the previous example in the following way.
    Let $\gfr = \gfr_0 + \hfr_0$ be a double Lie algebra and consider bases $\{X_1, \dots, X_n\} \subset \gfr_0$ and $\{Y_1, \dots, Y_m\} \subset \hfr_0$ of $\gfr_0$ and $\hfr_0$ respectively.
    \zcref[S]{th: double lie algs and lie bialgs of translation type} gives a Manin triple $(\Gfr = \gfr \ltimes_{\ad^{\flat}} \gfr^*_{\ab}, \gfr_0 \ltimes \hfr_{0,\ab}^*, \hfr_0 \ltimes \gfr_{0,\ab}^*)$ with
    the factorisable $r$-matrix
    \begin{align*}
        r = \sum_{i = 1}^n X_i \tensor X_i^* + \sum_{j= 1}^m Y_j^* \tensor Y_j.
    \end{align*}
    The double of \(\Gfr\) corresponds to the Manin triple
    \((\Dfr = \Gfr \times \Gfr, \diag(\Gfr), (r_+^{\sharp}, r_-^{\sharp})(\Gfr^*))\).
    If we now set $\dfr := \gfr \times \gfr$, we obtain the double Lie algebra with reductive complement
    \begin{align*}
        \dfr = \diag(\gfr) + (\hfr_0 \times \gfr_0) = \diag(\gfr) + (0 \times \gfr).
    \end{align*}
    This reduces to the previous example when $\hfr_0 = 0$.
\end{example}

Combining \zcref{prop: abelian approximation from lie--reductive doubles,cor: reductive complement from coboundary of translation type} leads to a general strategy to quantize triangular as well as some quasi-triangular Lie bialgebras of translation type.
We will further examine these situations in a following paper and construct a dual unitary \(2\)-cocycle in the triangular case and a quantum \(R\)-matrix in the quasi-triangular case.

\begin{corollary}\label{cor: abelian approximation for triangular TT}
    Consider a Lie algebra \(\Gfr = \gfr \ltimes_{\alpha} V\) and a triangular \(r\)-matrix of translation type \(r = r_{\gfr,V} + r_{V,\gfr} \in \Gfr \wedge \Gfr\), where \(r_{\gfr,V} \in \gfr \tensor V\), \(r_{V,\gfr} \in V \tensor \gfr\).
    Then \(\hfr = \{f + r^{\sharp}(f) : f \in V^*\} \subseteq \gfr \ltimes_{\alpha^{\flat}} V^*\) is a Lie subalgebra.
    If the double Lie algebra \(\gfr \ltimes_{\alpha^{\flat}}V^* = \gfr + \hfr\) integrates to a double Lie group \(G \ltimes V^* = G \bowtie H\), then the map
    \begin{align*}
        &R^{V^*} : H \to V^* \\
        &R^{V^*}(h) := f && \text{for all } h = gf \in H \subseteq G \ltimes V^*\!, g \in G, f \in V^*
    \end{align*}
    is an abelian approximation.
\end{corollary}

\begin{corollary}\label{cor: abelian approximation for quasi-triangular TT}
    Consider a Lie algebra \(\Gfr = \gfr \ltimes_{\alpha} V\) and a quasi-triangular \(r\)-matrix of translation type \(r = r_{\gfr,V} + r_{V,\gfr} \in \Gfr \wedge \Gfr\), where \(r_{\gfr,V} \in \gfr \tensor V\), \(r_{V,\gfr} \in V \tensor \gfr\).
    We write \(\dfr = \gfr + \hfr = \gfr + \pfr\) for the corresponding double Lie algebra with reductive complement.
    If the ideal \(\pfr \leq \dfr\) is the Lie algebra of an exponential Lie group \(P\) and there exists a double Lie group \(G \ltimes P = G \bowtie H\), then the map
    \begin{align*}
        \eta_H = \exp_{\pfr}^{-1} \circ R^P : H \to \pfr
    \end{align*}
    is an abelian approximation.
\end{corollary}

\section{Examples}\label{sec:examples}
We now discuss examples of abelian approximations and Poisson--Fourier transforms.
We first look at Lie groups $G$ equipped with the trivial Poisson bracket $\{\cdot, \cdot\}_G = 0$.
The quantum duality then recovers precisely the philosophy of Kirillov's orbit method, where we wish to construct the unitary dual of $G$ via quantization of $\gfr^*$.
This can be accomplished for exponential Lie groups, and we will reinterpret this result through the lens of the Poisson--Fourier transform.
We also discuss the Poisson double of \(G\), which is nothing else but the cotangent bundle \(T^* G = G \ltimes_{\Ad^*{\flat}} \gfr^*\) as we have seen in \zcref{example: cotangent bundle}.
Here we obtain matched abelian approximations from a \(G\)-equivariant measure class isomorphism \(\eta : G \to \gfr\), where we consider the conjugation action on \(G\).
The resulting quantization of \(T^*G\) provides a non-formal analogue of the Gutt \(\star\)-product \cite{guttExplicitProductCotangent1983}.

We then turn our attention to Lie groups with non-trivial Poisson structures.
We will explain how to realise the quantum duality principle for semi-direct product $G \ltimes_{\rho} V$ equipped with an essentially bijective $1$-cocycle $\gamma : G \to \hat V$, as well as for some two-step nilpotent groups.
Before considering the general case, we will investigate two key examples:
the quantum $ax+b$ group of Baaj and Skandalis \cite{skandalisDualityLocallyCompact1991} and the quantum Heisenberg groups of Szymczak and Zakrzewski \cite{SZ90QuantumHeisenberg}.
We will describe the Lie--Poisson structures we consider, the preferred Poisson dual group, the construction of the quantization, and the resulting Poisson--Fourier transform.
As we will see it is necessary in some cases to take pick a non-connected Poisson dual.
We want to argue with these examples that although the choice of the connected and simply connected dual is very convenient, it is not necessarily the best one when it comes to quantization.

We conclude by quickly mentioning the T-Poincaré models of Mercati \cite{mercatiTMinkowski_I_2024,mercatiTMinkowski_II_2024}.
These are a family of triangular Lie--Poisson structures of translation type of the Poincaré group \(\Pcal = \SO(1,3) \ltimes \R^4\).
We will come back to these examples a future article and present explicit matched pairs \(\SO(1,3) \bowtie H = \Pcal\) for each of them as in \zcref{cor: abelian approximation for triangular TT}.

\subsection{Exponential Lie groups and the orbit method}
Any Lie group $G$ can be turned into a Poisson--Lie group by equipping it with the trivial Poisson bracket $\{\cdot, \cdot\}_G = 0$.
The Poisson dual is then \(\gfr^*\) equipped with its linear Poisson structure.
There is no quantization of $\Lrm^{\infty}(G)$, and the quantum dual $W^*(G)$ models the tempered unitary dual of $G$.
Hence we wish to understand $W^*(G)$ as a quantization of $\gfr^*$.
This is very reminiscent of Kirillov's orbit method, which aims to construct the unitary dual $\hat G$ through quantization of $\gfr^*$ and its coadjoint orbits \cite{kirillovLecturesOrbitMethod2004}.

The programme of the orbit method has been successfully realised for exponential Lie groups, i.e., Lie groups $G$ whose exponential map $\exp : \gfr \to G$ is a diffeomorphism.
In this case there is an explicit homeomorphism $\gfr^* / G \to \hat G$ between the space of coadjoint obits and the unitary dual of $G$ \cite{arnalRepresentationStarExponentiels1990,leptinUnitaryRepresentationTheory1994}.
A key role is played by the pairing
\begin{align*}
    &\left(\cdot \mid \cdot\right) : \gfr^* \times G \to U(1) \\
    &\left(\xi \mid g\right) := e^{-i \la \xi, \log g\ra},
\end{align*}
where $\log : G \to \gfr$ denotes the inverse of $\exp$.
We can understand this pairing as the kernel of the Poisson--Fourier transform between $G$ and $\gfr^*$.
More precisely, set $H = V = 0$ to be the trivial group, and $\Gcal = G \ltimes V = G$, $\Hcal = H \ltimes \gfr^* = \gfr^*$.
Then $\eta_H = \id_0$ and $\eta_G = \log : G \to \gfr$ form matched abelian approximations for $G$ and $\gfr^*$ in the sense of \zcref{def:matchedAbApprox}.
The resulting Poisson--Fourier transform $\Fcal_G : \Lrm^2(G) \to \Lrm^2(\gfr^*)$ provides an isomorphism of \lcqgs
\[
    \Fcal_G : (W^*(G), \hat \Delta_G^{\alpha})\sphat\; = (\Lrm^{\infty}(G), \Delta_G) \cong (W^*(\gfr^*), \hat\Delta_{\gfr^*}^{\beta}).
\]
By applying duality, we then realise $(W^*(G), \hat \Delta_G)$ as a quantization of $(\Lrm^{\infty}(\gfr^*), \Delta_{\gfr^*})$.

\subsection{Cotangent bundles and quantum doubles}
Consider a Lie algebra \(\gfr\) equipped with the trivial cocommutator \(\gamma_{\gfr} = 0\).
As we discussed in \zcref{example: cotangent bundle}, the classical double of \((\gfr, \gamma_{\gfr})\) is given by the factorisable Lie algebra of translation type \(\Gfr = \gfr \ad^{\flat} \gfr^*_{\ab}\) with the \(r\)-matrix
\begin{align*}
    r = \sum_{i = 1}^n X_i \tensor X_i^*,
\end{align*}
where \(\{X_1, \dots, X_n\} \subset \gfr\) is any basis of \(\gfr\) and \(\{X_1^*, \dots, X_n^*\} \subset \gfr^*\) denotes the dual basis.
The corresponding double Lie algebra with reductive complement is
\begin{align*}
    \dfr = \gfr \ltimes_{\ad} \gfr = (\gfr \ltimes 0) + (0 \ltimes \gfr).
\end{align*}
Any Lie group \(G\) integrating \(\gfr\) then gives rise to the matched pair \(D = G \ltimes_c G = (G \ltimes e_G) \bowtie (e_G \ltimes G)\), where \(c_g(g_0) = g g_0 g^{-1}\) is the conjugation action.
An abelian approximation \(\eta : e_G \ltimes G \to \gfr\) is thus a measure-class isomorphism \(\eta : G \to \gfr\) which satisfies \(\eta(g_0 g_1 g_0^{-1}) = \Ad_{g_0} \eta(g_1)\).
Such maps are also known as \emph{generalised Cayley transforms} and have been studied before \cite{kostantGeneralizedCayleyMap2003,lemireCayleyGroups2006}.
The bicrossed product \(\hat G \bicrossl_{c} G\) coincides with the dual of quantum double \(\Dcal(G)\) of G.
\begin{proposition}\label{prop: bicrossed product and quantum double}
    Consider a \lcg \(G\) and write \(c : G \times G \to G : c_{g}(g_0) := g g_0 g^{-1}\) for the conjugation action.
    Then \(\hat G \bicrossl_{c} G = \Dcal(G)\sphat\), where \(\Dcal(G) = \Dcal(\Lrm^{\infty}(G))\) denotes the quantum double of \(G\) as described in \cite[Section~8]{baajDoubleCrossedProducts2005}.
\end{proposition}
\begin{proof}
    We know from \zcref{prop:dualityBicrossed} that the flip map \(\Sigma : \Lrm^2(G) \tensor \Lrm^2(G) \to \Lrm^2(G) \tensor \Lrm^2(G)\) induces an isomorphisms of \lcqgs
    \begin{align*}
        \Sigma : \left(\hat G \bicrossl_c G\right)\sphat \cong \hat G \bicrossl^c G.
    \end{align*}
    Therefore it is sufficient to show that \(\Sigma\) induces an isomorphism
    \begin{align*}
        \Sigma : \hat G \bicrossl^c G \cong \Dcal(G).
    \end{align*}
    This clearly holds at the level of von Neumann algebras, since \(\hat G \bicrossl^c G = W^*(G) \htensor \Lrm^{\infty}(G)\) and \(\Dcal(G) = \Lrm^{\infty}(G) \htensor W^*(G)\).
    We set \(\Zcal := W_G\) and \(Z = \Zcal (J_G \tensor \hat J_G) \Zcal (J_G \tensor \hat J_G)\).
    One easily computes that
    \begin{align*}
        \Zcal \xi(g_1, g_2) = \xi(g_1, g_1^{-1} g_2 g_1) \delta_G(g_1)^{1/2} && \text{for all } \xi \in \Lrm^2(G) \tensor \Lrm^2(G).
    \end{align*}
    We compute the matching \(m : \Lrm^{\infty}(G) \htensor W^*(G) \to \Lrm^{\infty}(G) \htensor W^*(G) : m(x) := Z x Z^*\).
    Clearly \(m(f \tensor 1) = f \tensor 1\) for all \(f \in \Lrm^{\infty}(G)\).
    Moreover,
    \begin{align*}
        m(1 \tensor \lambda_{g}) \xi(g_1, g_2)
            = (1 \tensor \lambda_g) \Zcal^* \xi(g_1, g_1^{-1} g_2 g_1) \delta_G(g_1)^{1/2}
            = \xi(g_1, g_1 g^{-1} g_1^{-1} g_2)
            = \lambda_{c_{g_1}(g)} \xi(g_1, g_2).
    \end{align*}
    The coproduct \(\Delta_m\) of \(\Dcal(G)\) is then determined by
    \begin{align*}
        &\Delta_m(f \tensor 1) = \Delta_G(f)_{13}; \\
        &\Delta_m(1 \tensor \lambda_g) \xi(g_1, g_2, g_3, g_4)
            = \xi(g_1, g_3 g^{-1} g_3^{-1} g_2, g_3, g^{-1} g_4).
    \end{align*}
    On the other hand using \cite[Propositions~2.4~and~2.5]{vaesExtensionsLocallyCompact2003}, we see that
    \begin{align*}
        &\Delta^c(1 \tensor f) = \Delta_G(f)_{24}; \\
        &\Delta^c(\lambda_g \tensor 1)\xi(g_1, g_2, g_3, g_4) = \xi(c_{g_4}(g^{-1}) g_1, g_2, g^{-1} g_3, g_4).
    \end{align*}
    Clearly \(\Sigma\) intertwines both of these coproducts, which concludes.
\end{proof}

Any abelian approximation \(\eta : G = (0 \ltimes G) \to \gfr\) also provides a matched abelian approximation \(\eta : (G \ltimes 0) = G \to \gfr\) with classical limit \(G \times \gfr^*\).
We obtain quantizations \((W^*(G \ltimes \gfr^*), \hat \Delta_{G \ltimes \gfr^*}^{\id_G})\) and \((W^*(G \times \gfr^*), \hat \Delta_{G \times \gfr^*}^c)\) together with a Poisson--Fourier transform
\begin{align*}
    \Fcal_{G \ltimes \gfr^*} : (W^*(G \ltimes \gfr^*), \hat \Delta_{G \ltimes \gfr^*}^{\id_G})\sphat \ \cong (W^*(G \times \gfr^*), \hat \Delta_{G \times \gfr^*}^c).
\end{align*}

\begin{example}
    If \(\gfr\) is the Lie algebra of an exponential Lie group \(G\), then \(\eta := \exp_{\gfr}^{-1} : G \to \gfr\) is an abelian approximation, see \zcref{cor: abelian approximation for quasi-triangular TT}.
\end{example}

\begin{example}
    We take \(G = \GL_n\) and identify \(\gfr = \gl_n\) with the space of \(n \times n\) matrices.
    Then the map
    \begin{align*}
        \eta_{\GL_n}(M) := M - 1_n
    \end{align*}
    is a \(G\)-equivariant measure class isomorphism, where \(1_n\) denotes the \(n \times n\) identity matrix.
\end{example}

The quantization \((W^*(G \times \gfr^*), \hat \Delta_{G \times \gfr^*}^c)\) is related with the famous Gutt \(\star\)-product on the cotangent bundle \cite{guttExplicitProductCotangent1983}.
Indeed let us define an associative product \(\star : \Acal(G \times \gfr^*) \times \Acal(G \times \gfr^*) \to \Acal(G \times \gfr^*)\) by
\begin{align*}
    (\vphi \star \psi)(x) := \la \vphi \tensor \psi, \hat \Delta^c_{G \times \gfr^*}(x) \ra && \text{for all } \vphi, \psi \in \Acal(G \times \gfr^*), x \in G \times \gfr^*.
\end{align*}
Then \(\star\) is a non-formal analogue of Gutt's product.

\subsection{The quantum \texorpdfstring{$ax+b$}{ax+b}-group of Baaj and Skandalis.}\label{sec:ax+b}
The quantum $ax+b$ group of Baaj and Skandalis was first presented in a talk at Oberwolfach using a bicrossed product \cite{skandalisDualityLocallyCompact1991}, and has been revisited several times.
In 2003, Vaes and Vainerman described the associated infinitesimal Hopf algebra and interpret it as a deformation of the $ax+b$ Lie algebra \cite{vaesExtensionsLocallyCompact2003}.
Ten years later, in 2013, Stachura constructed a dual unitary $2$-cocycle on the $ax+b$-group and conjectured the associated quantization to be isomorphic to the Baaj--Skandalis bicrossed product.
This was proven by Bieliavsky--Gayral--Neshveyev--Tuset in 2021 \cite{bieliavskyQuantizationSubgroupsAffine2021}.
See also \cite[Section 4]{baajNonSemiRegularQuantumGroups2003} for a discussion of this quantum group over an arbitrary locally compact ring.

We write \(G := \R^{\times}\) and consider the group $\Gcal = G \ltimes \R$, so that its group structure is given by 
\[
    (g_1, v_1) (g_2, v_2) = (g_1 g_2, v_1 g_2 + v_2)
\]
for all $(g_1, v_1), (g_2, v_2) \in G \ltimes \R$.
The Lie algebra $\gfr = \Lie(\Gcal)$ is generated by the elements $\partial_g$ and $\partial_v$, with bracket
\[
    [\partial_v, \partial_g] = \partial_v.
\]
For each \(\hbar \in \R\setminus\{0\}\), $\gfr$ admits a triangular non-degenerate $r$-matrix of translation type $r_{\hbar} = \hbar \partial_v \wedge \partial_g$.
We follow \zcref{cor: abelian approximation for triangular TT} to quantize \(\Gcal\).
Let \(\dfr := \gfr \ltimes \R^{\times}\).
It admits the basis \(\{\partial_g, \drm_v\}\) with bracket \([\partial_g, \drm_v] = \drm_v\).
The subspace \(\hfr = \{r_{\hbar}^{\sharp}(\xi) + \xi : \xi \in \R^*\}\) is therefore generated by the element \(\hbar \partial_g + \drm_v \in \dfr\).
We can integrate this structure to the Lie groups \(D = \R^* \rtimes G\) and \(H_\hbar = \{((g-1)\drm_v / \hbar, g) : g \in G\} \leq D\).
We clearly have a matched pair \(D = G \bowtie H_{\hbar}\), and by \zcref{cor: abelian approximation for triangular TT}, the map \(\Rrm^{\R^*}_{H,\hbar}((g-1)\drm_v / \hbar, g) = (1 - g^{-1})\drm_v / \hbar\) is an abelian approximation of \(H_{\hbar}\) with classical limit \(\Gcal\).

We construct a matched abelian approximation of \(G\) in a similar way by conjugating the matched pair \(D = G \bowtie H_{\hbar}\) with \((\drm_v/\hbar) \in D\).
Clearly, \((\drm_v/\hbar) H_{\hbar} (-\drm_v/\hbar) = G =: H\) and \(G_{\hbar} := (\drm_v/\hbar) G (-\drm_v/\hbar) =\{((1 - g)\drm_v / \hbar, g) : g \in \R^{\times}\}\).
In particular \(R^{\R^*}_{G_{\hbar}}((1 - g)\drm_v / \hbar, g) = (g^{-1} - 1)\drm_v / \hbar\) is an abelian approximation of \(G_{\hbar}\).
But \(D = G' \bowtie H' = G \bowtie H\) are isomorphic matched pairs, and therefore \(\eta_{G,\hbar}(g) = \frac{1}{\hbar}(g^{-1} - 1)\drm_v\) is an abelian approximation of \(G\) with classical limit \(\Hcal := H \ltimes \R\).
We identify \(H_{\hbar}\) with \(H := G\) using the conjugation and write \(\eta_{H,\hbar}(h) = \frac{1}{\hbar} (1-h^{-1}) \drm_v\).
We obtain quantizations $(W^*(\Gcal), \hat \Delta_{\Gcal,\hbar}^{\beta})$ of $(W^*(\Gcal), \hat \Delta_{\Gcal})$ and $(W^*(\Hcal), \hat \Delta_{\Hcal,\hbar}^{\alpha})$ of $(W^*(\Hcal), \hat \Delta_{\Hcal,\hbar})$ together with the Poisson--Fourier transform
\begin{align*}
    \Fcal_{\Gcal,\hbar} : (W^*(\Gcal), \hat \Delta_{\Gcal,\hbar}^{\beta})\sphat\ \cong (W^*(\Hcal), \hat \Delta_{\Hcal,\hbar}^{\alpha}).
\end{align*}
We note that the family \(\{(W^*(\Gcal), \hat \Delta_{\Gcal,\hbar}^{\beta})\}_{\hbar \neq 0}\) converges to \((W^*(\Gcal), \hat \Delta_{\Gcal})\) in the strong operator topology as \(\hbar \to 0\), see \cite[Corollary 3.27]{bieliavskyQuantizationSubgroupsAffine2021}.

\subsection{Symplectic Lie groups with Lagrangian normal subgroup}

The previous example generalises to a much larger class of \lcgs, as described in \cite{bieliavskyQuantizationSubgroupsAffine2021,bieliavskyQuantizationLocallyCompact2024}.
We consider the setting of \cite{bieliavskyQuantizationLocallyCompact2024} with trivial \(2\)-cocycle.
Let $\Gcal = G \ltimes_{\rho} V$ be a semi-direct product, where $V$ is abelian and $\rho : V \times G \to V$ is a right action by continuous automorphisms.
We equip $\hat V$ with the dual left action $\hat \rho$.
The following definition gives the main ingredient necessary for this construction:
\begin{definition}\label{def:essentially_bijective_1_cocycle}
    An \emph{essentially bijective $1$-cocycle} for the action $\hat \rho$ is a continuous $1$-cocycle $\gamma : G \to \hat V$ for $\hat \rho$ which is also a measure class isomorphism.
\end{definition}

Suppose from now on that we have an essentially bijective $1$-cocycle $\gamma : G \to \hat V$.
The dual group we consider is $\Hcal := \Gcal$.
For clarity, we set $H = G$, $W = V$, $\sigma = \rho$, and write $\Hcal = H \ltimes_{\sigma} W$.
Finally, set $D = \hat V \rtimes_{\hat \rho} G$, and $H_{\gamma} = \{(\gamma(h), h) : h \in H\}$.
Since $\gamma$ is a measure class isomorphism, we have a matched pair $G \bowtie H_{\gamma} = D$.
We compute the dressing transformations:
\begin{align*}
    &\mu(g,h) = (0, g)(\gamma(h^{-1}), h^{-1}) = (\hat\rho_g(\gamma(h^{-1})), g h^{-1}); \\
    &\mu_{\Sigma}(g,h) = (\gamma(h^{-1}), h^{-1})(0, g) = (\gamma(h^{-1}), h^{-1} g); \\
    &(\beta_h(g), \alpha_g(h)) = \mu_{\Sigma}^{-1} \mu (g,h)
        = \left(
            \gamma^{-1}(\hat\rho_g(\gamma(h^{-1})))^{-1} g h^{-1},
            \gamma^{-1}(\hat\rho_g(\gamma(h^{-1})))^{-1}
        \right).
\end{align*}
We further analyse $\beta_h(g)$:
\begin{align*}
    \gamma\left(\beta_h(g)^{-1}\right)
        &= \gamma\left(h g^{-1} \gamma^{-1}(\hat\rho_g(h^{-1}))\right) \\
        &= \gamma(h) + \hat\rho_h\left(\gamma(g^{-1}) + \hat\rho_{g^-1} \hat\rho_g(\gamma(h^{-1}))\right) \\
        &= \gamma(h) + \hat\rho_h(\gamma(g^{-1})) + \hat\rho_h(\gamma(h^{-1})) \\
        &= \hat\rho_h(\gamma(g^{-1})) + \gamma(h h^{-1}) \\
        &= \hat\rho_h(\gamma(g^{-1})).
\end{align*}
Hence, $\beta_h(g) = \gamma^{-1}(\hat\rho_h(\gamma(g^{-1})))^{-1} = \alpha_h(g)$ for almost all $g \in G, h \in H$.
If we now define $\eta_H : H \to \hat V : h \mapsto -\gamma(h^{-1})$ and $\eta_G : G \to \hat W : g \mapsto \gamma(g^{-1})$, then $(\eta_G, \eta_H)$ form matched abelian approximations for $(\Gcal, \Hcal)$.
Moreover since $\alpha_x(y) = \beta_x(y)$ for almost all $x, y \in G$, the Poisson--Fourier transform $\Fcal_{\Gcal}$ shows that $(W^*(\Gcal), \hat \Delta_{\Gcal}^{\beta})$ is self-dual:
\[
    \Fcal_{\Gcal} : \left(W^*(\Gcal), \hat \Delta_{\Gcal}^{\beta}\right)\sphat\; \xrightarrow{\cong}\; \left(W^*(\Hcal), \hat\Delta_{\Hcal}^{\alpha}\right) = \left(W^*(\Gcal), \hat \Delta^{\beta}_{\Gcal}\right).
\]

Besides the examples discussed in \cite[Section 2.4]{bieliavskyQuantizationSubgroupsAffine2021} and \cite[Examples 1.5, 1.6]{bieliavskyQuantizationLocallyCompact2024}, we can try to construct essentially bijective \(1\)-cocycles from symplectic Lie groups with a Lagrangian normal subgroup.
We already saw this strategy in action in the previous example maps in order to get the maps \(\eta_{G,\hbar}\) and \(\eta_{H,\hbar}\).
Consider a symplectic Lie algebra \((\Gfr, \omega)\) with a Lagrangian ideal \(L \leq \Gfr\).
Any such ideal is abelian \cite[Lemma~2.1.]{bauesSymplecticLieGroups2016}, so \(\Gfr\) is an extension
\begin{align*}
    V \to \Gfr \to \gfr := \Gfr / V
\end{align*}
We assume for simplicity that the above extension splits, although it is not required for the following arguments to go through.
It is well-known that \(r := \omega^{\sharp}(\omega)\) is a non-degenerate triangular \(r\)-matrix on \(\Gfr\), see for example \cite[Corollary~2.2.4]{chariGuideQuantumGroups1995}.
Hence by \zcref{cor: abelian approximation for triangular TT} the subspace
\begin{align*}
    \hfr = \{r^{\sharp}(\xi) + \xi : \xi \in \gfr^*\} = \{X + r^{\flat}(X) : X \in \gfr\} \leq \gfr \ltimes V^*
\end{align*}
is a Lie subalgebra, and \(j : \gfr \to \gfr \ltimes V^* : j(X) := X + r^{\sharp}(X)\) is a section of the quotient map \(\gfr \ltimes V^* \to \gfr\).
If \(j\) integrates to a Lie group homomorphism \(J : G \to G \ltimes V^* : J(g) := g \eta(g)\) such that \(G \ltimes V^* = G \bowtie J(G)\) forms a matched pair, then \(\eta : G \to V^*\) is an essentially bijective \(1\)-cocycle.

\subsection{The quantum Heisenberg groups of Szymczak and Zakrzewski}\label{sec:heisenberg}
We now turn our attention to the quantum Heisenberg groups constructed in \cite{SZ90QuantumHeisenberg}, which have been shown to come from a cocycle bicrossed product in \cite{zakrzewskiPoissonLieGroups1992}.

Let $V$ be a finite dimensional real vector space equipped with a symplectic form $\omega$, and choose $M \in \mathfrak{sp}(V, \omega)$.
That is, $M \in \End(V)$ such that $\omega(Mv, w) = \omega(Mw, v)$ for all $v, w \in V$.
We write $\pi_{\omega}$ for the constant Poisson structure on $V$ obtained from $\omega$.
We let $G = \R$, $\rho(v)_{g} := e^{g M}v$, and $\Gcal = G \ltimes_{\rho} V$.
Then $\pi_{\Gcal}(g,v) = g \pi_{\omega}$ defines a Lie--Poisson structure on $\Gcal$.
We identify $\gfr = \Lie(\Gcal)$ with $G \times V$ and set $\partial_g = (1, 0_V) \in \gfr$, so that the Lie bracket on $\gfr$ is given by
\begin{align*}
    [v, \partial_g]_{\gfr} = M v, && [v_1, v_2]_{\gfr} = 0, && \text{for all } v, v_1, v_2 \in V.
\end{align*}
The dual bracket $[-,-]_{\gfr^*}$ on $\gfr^*$ induced by $\pi_{\Gcal}$ is
\begin{align*}
    [\xi, \dd g]_{\gfr^*} = 0, && [\xi, \zeta]_{\gfr^*} = \omega(\omega^{\sharp}(\xi), \omega^{\sharp}(\zeta)) \cdot \dd g,
    && \text{for all } \xi, \zeta \in V^*,
\end{align*}
where $\dd g \in G^*$ is dual to $\partial_g$, and $\omega^{\sharp} : V^* \to V$ denotes the musical isomorphism.
We thus see that $(\gfr^*, [-,-]_{\gfr^*})$ is isomorphic to the Heisenberg Lie algebra.

More precisely, we set $H = V$, $W = \R$, $\psi(h_1, h_2) = \frac{1}{2} \omega(h_1, h_2)$, and consider the central extension $\Hcal = H \ltimes_{\id,\psi} W$.
Thus $\Hcal$ is the Heisenberg group on $(V,\omega)$.
We identify $\hfr = \Lie(\Hcal)$ with $H \times W$ and write $\partial_w = (0_H, 1) \in \hfr$, so that the Lie bracket on $\hfr$ is
\begin{align*}
    [h_1, h_2]_{\hfr} = \omega(h_1, h_2) \cdot \partial_w,
    && [h, \partial_w] = 0,
    && \text{for all } h, h_1, h_2 \in V.
\end{align*}
The map $\phi : \gfr^* \to \hfr$ determined by $\dd g \mapsto \partial_w,\; \xi \mapsto \omega^{\sharp}(\xi)$ for $\xi \in V^*$ is clearly an isomorphism of Lie algebras.
By transporting the cocommutator on $\gfr^*$ we obtain a Lie bialgebra structure on $\hfr$, whose corresponding Poisson bivector on $\Hcal$ is given by
\begin{align*}
    \pi_{\Hcal}(h, w) = \partial_w \wedge M h.
\end{align*}

For the simplicity of the formulas, we use $\omega^{\sharp}$ to identify $V^*$ with $V$.
Choose $\hbar \in \R\setminus\{0\}$ and let $\eta_G : G \to W^* = \R : \eta_G(g) = \frac{g}{\hbar}$ and $\eta_H : H \to V \cong V^* : \eta_H(h) = \frac{h}{\hbar}$.
An easy computation shows that $(\eta_G, \eta_H)$ form matched abelian approximations for $(\Gcal, \Hcal)$.
The only non-trivial action is $\alpha_g(h) = e^{-g M}h$, and $\psi$ induces the cocycle
\[
    \Vcal(h_1, h_2, g) = e^{\frac{-i}{2} \eta_G(g) \omega(h_1, h_2)} = e^{-\frac{ig}{2\hbar} \omega(h_1, h_2)}.
\]
We note that $J_{\Gcal} : \Gcal \to \gfr : (g, v) \mapsto (\eta_G(g), v)$ and $J_{\Hcal} : \Hcal \to \hfr : (h, w) \mapsto (\eta_H(h), w)$ are Poisson maps, where $\gfr$ and $\hfr$ are equipped with the linear Poisson bracket coming from the Lie bialgebra structure.

\subsection{Two-step Nilpotent groups}
We generalise the previous examples by considering a general two-step nilpotent \lcg $\Hcal$.
The construction presented here is inspired by those in \cite{rieffelSolvableQuantumGroups1989,vandaeleQuantumDeformationHeisenberg1991} and is related to some of the examples in \cite[Examples 8.26]{baajUnitairesMultiplicatifsDualite1993}.
Consider a central extension $\Hcal = H \ltimes_{\id,\psi} W$ with both $H$ and $W$ abelian, as well as a semi-direct product $\Gcal = G \ltimes_{\rho} \hat H$ for a right action by continuous automorphisms $\rho : \hat H \times G \to \hat H$.
We furthermore assume there exists a left action $\theta : G \times W \to W$ such that
\[
    \psi(\hat \rho_g(h_1), \hat \rho_g(h_2)) = \theta_g(\psi(h_1, h_2))
\]
for all $h_1, h_2 \in H, g \in G$.
Note that if $\psi$ is surjective, then $\theta$ is completely determined by $\hat \rho$.
We let $\hat \theta$ denote the dual right action of $G$ on $\hat W$.

In this setup, any continuous $1$-cocycle $\gamma : G \to \hat W$ gives rise to a cocycle $\Vcal(h_1, h_2, g) = \ppairing{\psi(h_1,h_2)}{\eta_G(g)}$ for the matched pair $G \bowtie H = \left(H \rtimes_{\hat \rho} G\right)$, and the identity map $\eta_H = \id_H$ is an abelian approximation of $H$.
If $\gamma$ is essentially bijective in the sense of \zcref{def:essentially_bijective_1_cocycle}, then $(\eta_G = \gamma, \eta_H = \id_H)$ form matched abelian approximations for $(\Gcal, \Hcal)$.

\begin{remark}
    If $\theta$ is the trivial action and $\gamma$ is essentially bijective, then $\gamma$ must be an isomorphism of \lcgs by the $1$-cocycle property.
    We then recover the situation discussed in \cite{vandaeleQuantumDeformationHeisenberg1991}.
    In particular the quantum Heisenberg groups of the previous section fit in this picture: $\psi$ is a symplectic form on $H$ and $\hat \rho$ acts by symplectic automorphisms.
\end{remark}

\begin{remark}
    Although formulated in a slightly different language, the conditions we use to obtain a cocycle bicrossed product were already described in \cite[Example 5.1.1]{vaesExtensionsLocallyCompact2003}.
\end{remark}

\begin{example}
    Consider $U := \R^n, H := U \times U^*$, and $\psi(u_1, \xi_1; u_2, \xi_2) = \la \xi_1, u_2 \ra - \la \xi_2, u_1 \ra \in \R =: W$.
    Let $G := \R^{\times}$, choose $\mu \in \R\setminus\{-1\}$, and set
    \[
        \hat \rho_g(u, \xi) = (g u, |g|^{\mu} \xi).
    \]
    Then
    \[
        \psi(\hat \rho_g(h_1), \hat \rho_g(h_2)) = g |g|^{\mu} \psi(h_1, h_2)
    \]
    for any $h_1, h_2 \in H, g \in G$, and
    \[
        \eta_G : G \to \R \cong \hat W : g \mapsto \frac{1}{\hbar (\mu + 1)} (g|g|^{\mu} - 1)
    \]
    is an essentially bijective $1$-cocycle for the induced action on $\hat W$.
\end{example}

\begin{example}
    Consider a Lie group $G$, and let $H = W = \gfr := \Lie(G)$ as an abelian groups.
    We write $\hat \rho = \Ad$ for the adjoint action.
    Then the Lie bracket $[-, -]$ on $\gfr$ induces a $2$-cocycle $\psi(h_1, h_2) := [h_1, h_2] \in W$ which is $\hat \rho$-equivariant, so that we may take $\theta = \alpha$.
    Then any element $\xi_0 \in \gfr^*$ induces a $1$-cocycle
    $\eta_G(g) = \rho(\xi_0)_g - \xi_0$, and we obtain a quantization of the cotangent bundle $\Trm^*G$.
    If we assume that $\xi_0$ has trivial stabiliser and an open and dense orbit, then $\eta_G$ is essentially bijective.
    This way, we obtain matched abelian approximation $(\eta_G, \eta_H = \id_{\gfr})$ for $(\Gcal = G \ltimes_{\rho} \gfr^*, \Hcal = \gfr \ltimes_{\id,\psi} \gfr)$.
\end{example}

\subsection{T-Poincaré groups}\label{ex: T-Poincare}
There is a long history in the physics literature to consider non-commutativity relations on spacetime.
The first model is due to Snyder following an idea of Heisenberg \cite{snyderQuantizedSpaceTime1947}, and many others have been considered since such as Moyal-type deformations \cite{akoforQUANTUMFIELDSGROENEWOLD2008} or \(\kappa\)-deformations \cite{lukierskiNewQuantumPoincare1992,lukierskiKappadeformationsHistoricalDevelopments2017}.
At the semi-classical level, a Poisson bracket on the Minkowski spacetime \(M = \R^4\) of the form
\begin{align*}
    \{x^{\mu}, x^{\nu}\}_M = \Theta^{\mu \nu}(x)
\end{align*}
will usually not be Poincaré invariant.
This can be corrected by equipping \(\Pcal\) with a non-trivial Poisson--Lie structure, so that the action map \(\Pcal \times M \to M\) becomes Poisson.
The Lie--Poisson structures on \(\Pcal\) were studied by Zakrzewski who came up with a partial classification \cite{zakrzewskiPoissonStructuresPoincare1997}, see also \cite{tolstoyTwistedQuantumDeformations2008} for minor corrections.
A special subclass of these Lie--Poisson structures are the T-Poincaré models of Mercati \cite{mercatiTMinkowski_I_2024,mercatiTMinkowski_II_2024}, which he considered for their desirable mathematical and physical properties.
They correspond exactly to the triangular structures induced by an \(r\)-matrix of the form
\begin{align*}
    r = \frac{-1}{2} \theta^{\mu \nu} P_{\mu} \wedge P_{\nu} + \frac{-1}{2} \omega^{\mu \rho \sigma} P_{\mu} \wedge M_{\rho \sigma} \in \pfr \wedge \pfr,
\end{align*}
where \(\pfr = \spn \{P_\mu, M_{\rho \sigma}\} = \Lie(\Pcal)\) is the Poincaré Lie algebra generated by infinitesimal translations \(\{P_\mu\}\) and Lorentz transformations \(\{M_{\rho \sigma}\}\).
These are clearly triangular \(r\)-matrices of translation type by \zcref{prop: coboundary structures of translation type}, and we may hope to quantize them using \zcref{cor: abelian approximation for triangular TT} when \(\theta = 0\).
This project can be fully carried our in all of the T-Poincaré models, and we will present explicit solutions in an upcoming article.

\section*{Acknowledgements}
I am very grateful for the continuous support and guidance of my adviser Pierre Bieliavsky, and our many insightful conversations.
I would also like to thank S. Neshveyev and L. Tuset for their invitation while I was conducting this research, and V. Gayral for his helpful answers to my questions and his many recommendations.
A. Massar is a FRIA grantee of the Fonds de la Recherche
Scientifique -- FNRS.

\section*{AI usage disclosure}
The author believes that large AI companies pose real threats to our general well-being, including but not limited to privacy, cybersecurity, labour, the environment, and democratic processes at the local, national, and international level.
It is therefore their commitment to exclusively use these companies' models to find references to the existing mathematics literature.
In particular, no LLM was used while conducting the research in this article.

\printbibliography

\end{document}